\pdfoutput=1
\RequirePackage{ifpdf}
\ifpdf % We are running pdfTeX in pdf mode
\documentclass[pdftex]{sigma}
\else
\documentclass{sigma}
\fi

\usepackage{xspace}
\usepackage{bbm}

\numberwithin{equation}{section}

\newtheorem{Theorem}{Theorem}[section]
\newtheorem{Corollary}[Theorem]{Corollary}
\newtheorem{Lemma}[Theorem]{Lemma}
\newtheorem{Proposition}[Theorem]{Proposition}
{ \theoremstyle{definition}
\newtheorem{Definition}[Theorem]{Definition}

\newtheorem{Remark}[Theorem]{Remark} }

\DeclareMathOperator{\RE}{Re}
\DeclareMathOperator{\IM}{Im}

\newcommand{\mi}{\mathbbm{i}\xspace}
\newcommand{\bbC}{\mathbb{C}}
\newcommand{\bbE}{\mathbb{E}}
\newcommand{\bbH}{\mathbb{H}}
\newcommand{\bbN}{\mathbb{N}}
\newcommand{\bbR}{\mathbb{R}}
\newcommand{\bbS}{\mathbb{S}}
\newcommand{\bbZ}{\mathbb{Z}}
\newcommand{\R}{\mathbb{R}}
\newcommand{\C}{\mathbb{C}}
\newcommand{\one}{\mathbbm{1}}
\newcommand{\CPone}{\bbC\mathrm{P}^1}

\begin{document}
\allowdisplaybreaks

\newcommand{\arXivNumber}{1801.07032}

\renewcommand{\PaperNumber}{011}

\FirstPageHeading

\ShortArticleName{On Closed Finite Gap Curves in Spaceforms I}

\ArticleName{On Closed Finite Gap Curves in Spaceforms I}

\Author{Sebastian KLEIN~$^\dag$ and Martin KILIAN~$^\ddag$}

\AuthorNameForHeading{S.~Klein and M.~Kilian}

\Address{$^\dag$~Lehrstuhl f\"ur Mathematik~III, Universit\"at Mannheim, B 6, 28--29, 68131 Mannheim, Germany}
\EmailD{\href{mailto:s.klein@math.uni-mannheim.de}{s.klein@math.uni-mannheim.de}}

\Address{$^\ddag$~Department of Mathematics, University College Cork, Ireland}
\EmailD{\href{mailto:m.kilian@ucc.ie}{m.kilian@ucc.ie}}

\ArticleDates{Received June 14, 2019, in final form February 28, 2020; Published online March 04, 2020}

\Abstract{We show that the spaces of closed finite gap curves in ${\mathbb R}^3$ and ${\mathbb S}^3$ are dense with respect to the Sobolev $W^{2,2}$-norm in the spaces of closed curves in ${\mathbb R}^3$ respectively ${\mathbb S}^3$.}


\Classification{53A04; 37K10; 30D15; 46E35; 22E46}

\section{Introduction}

The shape of a curve $\gamma=\gamma(t)$ in a 3-dimensional space form is determined by two real-valued functions, the (geodesic) curvature $\kappa(t)$ and the torsion $\tau(t)$. Following an idea due to Hasimoto~\cite{hasimoto1972}, the information of these two functions is merged in a single complex-valued function, the complex curvature (Hasimoto map) $q(t) = \kappa(t)\exp\big(\mi \int_0^t \tau(s)\mathrm{d}s\big)$. The problem of reconstructing the extended frame $F=F(t,\lambda)$ of $\gamma$ from its complex curvature leads to a linear differential equation $\big(\tfrac{\mathrm{d}\ }{\mathrm{d}t} - \alpha\big)F=0$, where the operator $\tfrac{\mathrm{d}\ }{\mathrm{d}t} - \alpha$ is the differential operator of the 2-dimensional self-focusing nonlinear Schr\"odinger equation (NLS) with the potential $q$.

The NLS hierarchy, this is the hierarchy of flows induced by the NLS operator, is a completely integrable system (Zakharov and Shabat \cite{MR0406174})
that is related to the vortex filament hierarchy via the Hasimoto map \cite{hasimoto1972}.
The NLS integrable system has been the subject of very intensive research, it is probably the second-best researched infinite-dimensional completely integrable system (after the system of
the Korteweg--de Vries (KdV) equation). For an extensive overview of this research with many references, see Grinevich and Santini \cite[Introduction]{Grinevich/Santini}.
For curves in 2-dimensional spaceforms, the NLS flows reduce to the modified KdV (mKdV) hierarchy, and there is a wealth of literature also on this topic, see for example \cite{doliwa-santini1994,drinfeld-sokolov1981, goldstein-petrich1991, langer1999, langer-perline1998, langer-singer1984b} and the references therein.

In the present situation, \emph{closed} curves $\gamma$ are of particular interest. It should be noted that even if $\gamma$ is closed, and hence the functions $\kappa$ and $\tau$ are periodic, the complex curvature $q$ and therefore also the differential operator $\tfrac{\mathrm{d}\ }{\mathrm{d}t} - \alpha$ is generally only quasi-periodic. Conversely, even if~$q$ is periodic (this is called ``intrinsic periodicity''), the corresponding curve $\gamma$ is not necessarily closed. This is the case only if additional ``closing conditions'' (or conditions of ``extrinsic periodicity'') are satisfied. These closing conditions depend on the specific 3-dimensional space form we are considering, and they can be expressed in terms of the monodromy of the extended frame $F$.

Of particular interest in any integrable system are solutions which are stationary under all but finitely many of the flows in the hierarchy. Such solutions are called finite gap solutions. Such solutions are of interest in particular because they are solutions of an ordinary differential equation (of finite order). We thus say that a curve in a 3-dimensional space form is a finite gap curve if its complex curvature is stationary under all but finitely many of the flows in the NLS hierarchy. Grinevich and Schmidt \cite{grinevich-schmidt1995,grinevich-schmidt2000} introduced a deformation of finite gap curves that preserves both intrinsic and extrinsic periodicity. Callini and Ivey \cite{calini-ivey2007} consider the cabelling configurations of multiply-wrapped circles under these deformations.

One generally expects for a completely integrable system that the finite gap solutions are dense in the set of all solutions. The first result of this kind was shown by Mar\v{c}enko \cite{MR0374715} for the KdV system. Grinevich \cite{grinevich2001} proved that closed finite gap curves in $\mathbb{R}^3$ are dense in the set of all closed curves in $\mathbb{R}^3$. Our aim is to generalize his result and show that it also holds for curves in the 2-dimensional space forms $\bbR^2$, $\bbS^2$, $\bbH^2$ as well as in the other 3-dimensional space forms $\bbS^3$ and $\bbH^3$, and to present a new method of proof. In this first part we prove that closed finite gap curves in $\bbR^3$ and $\bbS^3$ are dense in the set of all closed curves in $\mathbb{R}^3$ respectively~$\bbS^3$. In a subsequent second part \cite{dense2} we deal with the remaining cases.

To prove the result for $\bbR^3$, Grinevich \cite{grinevich2001} uses the isoperiodic deformation of Grinevich and Schmidt \cite{grinevich-schmidt1995} and the Dubrovin equations
for reconstructing an NLS potential from the associated spectral data. Our proof is based on a different approach, utilising so-called perturbed Fourier coefficients
to characterise the NLS potentials, and applying asymptotic methods to these perturbed Fourier coefficients.
As far as we know, this kind of perturbed Fourier coefficients have not yet appeared in the study of the NLS integrable system before, however a related (but not identical) concept is that of
Birkhoff coordinates, which are discussed for example in \cite{Kappeler/Lohrmann/Topalov/Zung}.
We hope that the methods we here describe for finite gap curves can also serve as a ``model'' for the investigation of finite gap solutions in other integrable systems.

More specifically, the strategy for our proof is as follows: As possible complex curvature functions (potentials) we consider $L^2$-functions on an interval $[0,T]$; the corresponding curves are in the Sobolev space $W^{2,2}\big([0,T];\bbE^3\big)$. Given such a potential $q$, we use the ``monodromy'' $M(\lambda) := F(T,\lambda)$ of the corresponding initial value problem $\big(\tfrac{\mathrm{d}\ }{\mathrm{d}t} - \alpha\big)F=0$, $F(0,\lambda)=\one$, to construct a sequence $(z_k)_{k\in \bbZ}$ of ``perturbed Fourier coefficients'' of $q$ (see Definition~\ref{def:perturbedfourier}). Note that we carry out this construction even if $q$ is not periodic (in which case $M(\lambda)$ is not actually a~monodromy of the initial value problem). The benefit of these perturbed Fourier coefficients for the present paper lies in the fact that $q$ is finite gap if and only if all but finitely many of the~$z_k$ are zero, see Proposition~\ref{prop:fourier-finitegap}. The study of $q$ by means of its perturbed Fourier coefficients is made feasible by the second important instrument in this paper: an estimate of the difference between the monodromy $M(\lambda)$ of the given $q$ and the monodromy of the ``vacuum'' $q \equiv 0$, quantifying the idea that this difference becomes small for $|\lambda|$ large, see Theorem~\ref{thm:asymp}.
This asympotic estimate is obtained by methods that were developed by the first author in his thesis of habilitation \cite{klein-habil} for the monodromy of the sinh-Gordon equation.
Among the consequences of this estimate is that the sequence $(z_k)$ is $\ell^2$-summable, and that it is asymptotically close to the usual Fourier coefficients of $q$ (up to a constant factor); the latter fact justifies the name ``perturbed Fourier coefficients''.

\looseness=1 Drawing upon the asymptotic analysis of the monodromy, we show a partial analogue for perturbed Fourier coefficients to the well-known result that there is a 1--1 correspondence between $L^2$-functions and their (usual) Fourier coefficients: There exist finite-codimensional hyperplanes in $L^2([0,T])$ through the given $q$ and in $\ell^2(\bbZ)$ through $(z_k)$, so that the map associating
to each potential its perturbed Fourier coefficients is a local diffeomorphism between these hyperplanes (Proposition~\ref{prop:Phi-diffeo}). It follows from this result that the finite gap potentials
are dense in~$L^2([0,T])$ (Corollary~\ref{cor:dense}). Finite gap potentials are smooth and intrinsically periodic (Proposition~\ref{prop:finitegap-periodic}).

To obtain finite gap potentials which correspond to closed curves, we additionally need to satisfy the closing conditions besides the finite gap condition. In Section~\ref{Se:closing} we show
how to do so for curves in $\bbS^3$ and $\bbR^3$. The final result of this paper (Theorem~\ref{T:curvedense}) shows that the closed finite gap curves of length $T$ are dense in the Sobolev space of all closed $W^{2,2}$-curves of length~$T$ (in $\bbS^3$ or $\bbR^3$).

\section{Extended frames of curves}\label{Se:frames}

We study curves in the space forms $\bbR^3$ and $\bbS^3$. We identify Euclidean three space $\bbR^3$ with the matrix Lie algebra $\mathfrak{su}_2$ of skew-hermitian trace-free $2\times 2$
matrices. Under this identification we have
\[
	\| X \| = \sqrt{\det X},\qquad \langle X, Y \rangle = -\tfrac{1}{2} \mathrm{tr }(XY) \qquad \mbox{and} \qquad X \times Y = \tfrac{1}{2} [X,Y].
\]
The double cover of the isometry group under this
identification is $\mathrm{SU}_2 \ltimes \mathfrak{su}_2$.

We identify the three-sphere $\bbS^3 \subset \bbR^4$ with
$\bbS^3 \cong (\mathrm{SU}_2 \times \mathrm{SU}_2) / \mathrm{D}$, where D is the diagonal. The double cover of the isometry group
$\mathrm{SO}_4$ is $\mathrm{SU}_2 \times \mathrm{SU}_2$ via the action
$X \mapsto FXG^{-1}$. Under this identification, $\langle \cdot , \cdot \rangle$ is
the bilinear extension of the Euclidean inner product of
$\bbR^4$ to $\bbC^4$.

We fix a basis of $\mathfrak{sl}_2(\bbC)$ by
\begin{gather*} %\label{eq:epsilons}
 \varepsilon_- = \begin{pmatrix}
 0 & 0 \\ -1 & 0 \end{pmatrix},\qquad
 \varepsilon_+ = \begin{pmatrix}
 0 & 1 \\ 0 & 0 \end{pmatrix} \qquad \mbox{and}\qquad
 \varepsilon = \begin{pmatrix}
 \mi & 0 \\ 0 & -\mi \end{pmatrix}.
\end{gather*}
Each space form is endowed with a metric of an ambient
vector space, and we will denote by~$\langle \cdot , \cdot \rangle$
also the bilinear extension of the Ad-invariant inner
product of $\mathfrak{su}_2$ to $\mathfrak{su}_2^{\mbox{\tiny{$\bbC$}}} = \mathfrak{sl}_2(\bbC)$ such that
$\langle \varepsilon,\varepsilon \rangle = -\tfrac{1}{2} \operatorname{tr}\varepsilon^2 = 1$. We further have
\begin{gather*}
\langle \varepsilon_-,\varepsilon_- \rangle =
 \langle \varepsilon_+,\varepsilon_+ \rangle = 0,\qquad
 \varepsilon_-^* = - \varepsilon_+,\qquad \varepsilon^* = -\varepsilon, \nonumber\\
[\varepsilon_- ,\varepsilon ] = 2\mi \varepsilon_-,\qquad
 [\varepsilon, \varepsilon_+] = 2\mi \varepsilon_+ \qquad \mbox{and} \qquad
 [\varepsilon_+, \varepsilon_-]= \mi\varepsilon. %\label{eq:commutators}
 \end{gather*}

For a unit-speed curve $\gamma(t)$ in $\R^3$ or in $\bbS^3$, its (geodesic) curvature $\kappa(t)$ and its torsion $\tau(t)$ are defined in terms of the Frenet frame of $\gamma$.
Hasimoto~\cite{hasimoto1972} was the first to realise that if one replaces the normal components of the Frenet frame of $\gamma$ by a complex normal vector field~$Z_\gamma$
that is parallel with respect to the canonical connection of the complexified normal bundle of $\gamma$, one obtains an analogue to the Frenet--Serret equations for the Frenet frame.
In these equations, the two functions $\kappa$ and $\tau$ are replaced by the complex-valued function
\[
q(t) = \kappa(t) \exp\left( \mi \int_0^t \tau(s){\rm d}s \right) ,
\]
which is called the \emph{complex curvature} (or \emph{Hasimoto curvature}) of $\gamma$. Because of the existence of the analogue of the Frenet--Serret equations for $Z_\gamma$, we can expect that there is also an analogue of the Frenet--Serret theorem of curve theory in the new set-up,
meaning that a unit-speed curve~$\gamma$ is uniquely determined by its complex curvature up to a rigid motion of the ambient space. That this is indeed the case is shown in the following Lemma~\ref{th:sym}.
The complex curvature will be the main instrument with which we study curves in the 3-dimensional space forms.

For demonstration, we describe the construction of the complex curvature in $\R^3$, compare the reference \cite[Section~2]{hasimoto1972}. An analogous calculation applies for $\bbS^3$.
So consider a unit-speed curve~$\gamma(t)$ in~$\R^3$ and its Frenet frame $(T_\gamma,N_\gamma,B_\gamma)$. Here $T_\gamma=\gamma'$ and the normal fields~$N_\gamma$ and~$B_\gamma$,
the curvature $\kappa$ and the torsion $\tau$ are characterised by the Frenet--Serret equations
\begin{gather} \label{eq:frenet-serret}
 T_\gamma' = \kappa N_\gamma ,\qquad N_\gamma' = -\kappa T_\gamma + \tau B_\gamma \qquad\text{and}\qquad B_\gamma' = -\tau N_\gamma .
\end{gather}
Note that these equations show that the normal fields $N_\gamma$ and $B_\gamma$ are not parallel in the normal bundle $\perp\!\!(\gamma)$ of $\gamma$
equipped with the canonical connection (unless $\tau=0$). However, it also follows from equations~\eqref{eq:frenet-serret} that the section
\begin{gather*}
 Z_\gamma(t) = \bigr( N_\gamma(t) + \mi B_\gamma(t)\bigr) \exp\left( \mi \int_0^t \tau(s){\rm d}s \right)
\end{gather*}
of $\perp\!\!(\gamma) \otimes \C$ satisfies
\begin{gather} \label{eq:frenet-serret-complex}
 T_\gamma' = \RE(q \bar{Z}_\gamma) = \tfrac12\big(\bar{q}Z_\gamma + q\bar{Z}_\gamma\big) \qquad\text{and}\qquad Z_\gamma' = -q T_\gamma
\end{gather}
with the complex curvature $q(t) = \kappa(t) {\rm e}^{\mi \int_0^t \tau(s){\rm d}s}$.
Equations~\eqref{eq:frenet-serret-complex} are the analogue to the Frenet--Serret equations in our new set-up, and it follows from the equation for $Z_\gamma'$ that
$Z_\gamma$ has the advantage over $(N_\gamma,B_\gamma)$ that it is parallel with respect to the canonical connection of $\perp\!\!(\gamma) \otimes \C$.

We now return to unit-speed curves $\gamma$ in either $\R^3$ or $\bbS^3$.
If $\gamma$ is in the Sobolev space $W^{k,2}$ (meaning that $\gamma$ is $k$-times differentiable in the Sobolev sense, where the final derivative is square-integrable) for $k\geq 2$,
then its complex curvature is in~$W^{k-2,2}$. This is true even for $k=2$, where $\tau$ can only be regarded as a distribution in the negative Sobolev space $W^{-1,2}$ so that
$\int_0^t \tau(s){\rm d}s \in L^2$ holds, because even then ${\rm e}^{\mi \int_0^t \tau(s){\rm d}s}$ is bounded as $\tau$ is real-valued.

Also note that by the Sobolev embedding theorem, $W^{k,2}$ is contained in $C^{k-1}$ for any $k\geq 2$. In particular it follows that the curves~$\gamma$ constructed in the following lemma,
which is a variant of the Frenet--Serret theorem in our set-up, are at least continuously differentiable (in the classical sense) and thus have a continuous tangent map~$\gamma'$, regardless
of the value of $k\geq 2$. Therefore the evaluation of~$\gamma$ or~$\gamma'$ at specific points $t\in [0,T]$ is well-defined.

\begin{Lemma} \label{th:sym}
Let $k \in \{2,3,\dots,\infty\}$, $T>0$, $q \in W^{k-2,2}([0,T],\C)$ and
\begin{gather} \label{eq:alpha_q}
	\alpha^q = \tfrac{1}{2} \left( \lambda \varepsilon + q\varepsilon_+ + \bar{q}\varepsilon_- \right).
\end{gather}
Let $F = F(t,\lambda)$ solve
\begin{gather}\label{eq:frame1}
 \tfrac{{\rm d} \ }{{\rm d}t} F = F \alpha^q ,\qquad F(0,\lambda) = \mathbbm{1} \qquad \mbox{for all} \quad \lambda \in \mathbb{C}.
\end{gather}
Then $F$ depends holomorphically on $\lambda\in\bbC$ and we have $F(\cdot,\lambda) \in W^{k-1,2}\big([0,T],\bbC^{2\times 2}\big)$ for any $\lambda\in\C$. Moreover:
\begin{enumerate}\itemsep=0pt
\item[{\rm{(i)}}] $\gamma(t) = F(t,1)F(t,-1)^{-1}$ is a unit-speed curve in $\bbS^3$ with complex curvature $q$, with $\gamma \in W^{k,2}\big([0,T],\bbS^3\big)$ and with $\gamma'(t) = F(t,1)\varepsilon F(t,-1)^{-1}$.
\item[{\rm{(ii)}}] Write $F' = \tfrac{\partial F}{\partial \lambda}$. Then $\gamma(t)= 2 F'(t,0) F(t,0)^{-1}$ is a unit-speed curve in $\bbR^3$ with complex curvature~$q$, with $\gamma \in W^{k,2}\big([0,T],\bbR^3\big)$ and with $\gamma'(t) = F'(t,0)\varepsilon F(t,0)^{-1}$.
\end{enumerate}
\end{Lemma}

\begin{proof}
By definition, we have $\tfrac{\rm d}{{\rm d}t}F(\cdot,\lambda)\in W^{k-2,2}$ and therefore $F(\cdot,\lambda)\in W^{k-1,2}$.
Because $\alpha^q$ depends holomorphically on $\lambda$, $F$ also depends holomorphically on~$\lambda$.
By differentiating $\gamma(t)$ in~(i) and (ii), one can check
that $\gamma$ is a unit speed curve in $\bbS^3$ and in~$\R^3$ with complex curvature~$q$. $\gamma$ is a priori only in $W^{k-1,2}$, however by differentiating
the definition of $\gamma(t)$ with respect to $t$, one also obtains the formulas $\gamma'(t) = F(t,1)\varepsilon F(t,-1)^{-1}$ (for (i)) and $\gamma'(t) = F'(t,0)\varepsilon F(t,0)^{-1}$ (for (ii)),
which shows that $\gamma'(t) \in W^{k-1,2}$ and hence $\gamma \in W^{k,2}$ holds.
\end{proof}

In the situation of Lemma~\ref{th:sym}(i), where we consider curves in $\bbS^3$, we can replace the Sym points $\lambda_1=1$ and $\lambda_2=-1$ by any $\lambda_1 \neq \lambda_2 \in \bbR$,
giving the family of curves $\gamma_{\lambda_1,\lambda_2}(t) = F(t,\lambda_1)F(t,\lambda_2)^{-1}$ in $\bbS^3$. We have that
$\gamma'_{\lambda_1,\lambda_2}(t)=\tfrac12 (\lambda_1-\lambda_2)F(t,\lambda_1)\varepsilon F(t,\lambda_2)^{-1}$, therefore $\gamma_{\lambda_1,\lambda_2}$ has constant speed
$\tfrac12 |\lambda_1-\lambda_2|$ (unit speed for $|\lambda_1-\lambda_2|=2$).
Similarly in Lemma~\ref{th:sym}(ii) (for curves in $\bbR^3$) one can replace the Sym point $\lambda_0=0$ by any $\lambda_0 \in \bbR$, giving the family of curves
$\gamma_{\lambda_0}(t) = 2F'(t,\lambda_0) F(t,\lambda_0)^{-1}$. Again $\gamma_{\lambda_0}$ is a curve in $\bbR^3$ with
$\gamma'_{\lambda_0}(t) = F(t,\lambda_0)\varepsilon F(t,\lambda_0)^{-1}$ and hence, unit speed. In both cases, the family of these curves is called the
\emph{associated family} of curves corresponding to the original $\gamma$ or to $q$. We will see below that the translation of the Sym point or points along the real axis corresponds
to the changing of the total torsion of the associated curve.

The map $F\colon \mathbb{R} \times \mathbb{C} \to \mathrm{SL}_2(\bbC)$ is called an {\emph{extended frame}} of the curve. For any $k\geq 2$, the map $F(\cdot,\lambda)$ is at least continuous
by the Sobolev embedding theorem, and therefore the evaluation expression $F(t,\lambda)$ is well-defined for any $t\in [0,T]$ and $\lambda\in\bbC$.
Recall that $F(t,\lambda)$ depends holomorphically on $\lambda$. The extended frame $F$
satisfies the {\emph{reality condition}}
\begin{gather} \label{eq:reality}
	\overline{F\big(t,\bar\lambda\big)}^t = F^{-1}(t,\lambda) \qquad \mbox{for all}\quad \lambda \in \bbC \quad\mbox{and all} \quad t \in \mathbb{R}.
\end{gather}
In particular for $\lambda \in \bbR$ it takes values in $\mathrm{SU}_2$.

Now suppose $q\in L^2([0,T],\C)$ is the complex curvature of a closed, unit speed curve $\gamma$ in $\bbS^3$ or $\R^3$ of length $T$. Then $q$ is in general not periodic, but only quasi-periodic, meaning that
\begin{gather}\label{eq:q-quasiperiodic}
q(t+T)=\exp(\mi\theta T) q(t) \qquad \text{for almost all $t\in \R$,}
\end{gather}
with a unique number $\theta \in \R$. If $\gamma$ is three times differentiable, so that its torsion function $\tau$ is well-defined, then $\theta = \tfrac{1}{T} \int_0^T \tau(s) \mathrm{d}s \in \R$ is the average torsion of the curve $\gamma$, that is the total torsion of $\gamma$ divided by its length.
We will take the liberty of calling the number $\theta T$ the \emph{total torsion} of~$\gamma$ even where~$\gamma$ is only twice differentiable (and therefore
the torsion function of $\gamma$ is not in general well-defined).

In order to apply spectral theory to $\gamma$, we need a periodic potential. A periodic potential $\widetilde{q}$ can be obtained from $q$ by the regauging we now describe, compare for example
\cite[p.~11]{grinevich-schmidt-sfb}. We will see that $\widetilde{q}$ contains all the information which determines $\gamma$ \emph{except} for its total torsion, in other words
$\gamma$ can be reconstructed uniquely from $(\widetilde{q},\theta)$ (up to a rigid motion of the ambient space).

We regauge from the right by the matrix
\begin{gather*}
g(t) := \begin{pmatrix} \exp(\mi \theta t/2) & 0 \\ 0 & \exp(-\mi \theta t/2) \end{pmatrix} .
\end{gather*}
Note that $g$ depends on $q$ only via the value of $\theta$, it is independent of $\lambda$, and $g(0)=\one$ holds. With the regauged data
\begin{gather*}
 \widetilde{\alpha} := \alpha^q.g = g^{-1}\alpha^qg + g^{-1} \left( \tfrac{\mathrm{d}\ }{\mathrm{d}t} g \right) \qquad\text{and}\qquad \widetilde{F}(t,\lambda) = g^{-1}(0)F(t,\lambda)g(t) = F(t,\lambda)g(t) ,
\end{gather*}
$\widetilde{F}$ solves the analogous initial value problem to \eqref{eq:frame1}
\begin{gather*}%\label{eq:frame-periodic}
 \tfrac{{\rm d} \ }{{\rm d}t} \widetilde{F} = \widetilde{F} \widetilde{\alpha} ,\qquad \widetilde{F}(0,\lambda) = \mathbbm{1} \qquad \mbox{for all} \quad \lambda \in \mathbb{C}.
\end{gather*}
An explicit calculation shows that $\widetilde{\alpha}$ has the form of \eqref{eq:alpha_q}, meaning
\begin{gather*}%\label{eq:alpha-periodic}
 \widetilde{\alpha} = \tfrac{1}{2} \big( \widetilde{\lambda} \varepsilon + \widetilde{q}\varepsilon_+ + \overline{\widetilde{q}}\varepsilon_- \big) = \alpha^{\widetilde{q}}_{\widetilde{\lambda}}
\end{gather*}
with
\begin{gather*}
 \widetilde{q}(t) := \exp(-\mi \theta t)q(t)\quad \text{(for almost all $t\in [0,T]$)} \qquad\text{and}\qquad \widetilde{\lambda} = \lambda+\theta .
\end{gather*}
Note that $\widetilde{q}$ is periodic due to equation~\eqref{eq:q-quasiperiodic}, and that the translation of the spectral parame\-ter~$\lambda$ involved in this transformation is a measure of the total torsion~$\theta T$. Moreover, it follows from Lemma~\ref{th:sym} that the original curve $\gamma$ can be reconstructed from the data $(\widetilde{q},\theta)$: Let $F^{\widetilde{q}}$ be the extended frame defined by $\widetilde{q}$, that is the solution of the initial value problem
\begin{gather*}
 \tfrac{{\rm d} \ }{{\rm d}t} F^{\widetilde{q}} = F^{\widetilde{q}} \alpha^{\widetilde{q}} ,\qquad F^{\widetilde{q}}(0,\widetilde{\lambda}) = \mathbbm{1} \qquad \mbox{for all} \quad \widetilde{\lambda} \in \mathbb{C}.
\end{gather*}
We then have $F^{\widetilde{q}}\big(t,\widetilde{\lambda}\big) = \widetilde{F}\big(t,\widetilde{\lambda}-\theta\big) = F^q\big(t,{\widetilde{\lambda}-\theta}\big)g(t)$, and therefore in $\bbS^3$
\begin{gather}\label{eq:gamma-S3}
\gamma(t) = F^{\widetilde{q}}(t,1+\theta) F^{\widetilde{q}}(t,-1+\theta)^{-1} \qquad\!\text{and}\!\qquad \gamma'(t) = F^{\widetilde{q}}(t,1+\theta) \varepsilon F^{\widetilde{q}}(t,-1+\theta)^{-1} ,\!\!\!
\end{gather}
and in $\R^3$
\begin{gather}\label{eq:gamma-R3}
\gamma(t) = 2\left( \left. \tfrac{\partial F^{\widetilde{q}}(t,\widetilde{\lambda})}{\partial \widetilde{\lambda}}\right|_{\widetilde{\lambda}=\theta} \right) F^{\widetilde{q}}(t,\theta)^{-1}
\qquad\text{and}\qquad \gamma'(t) = \left( \left. \tfrac{\partial F^{\widetilde{q}}(t,\widetilde{\lambda})}{\partial \widetilde{\lambda}}\right|_{\widetilde{\lambda}=\theta} \right) \varepsilon F^{\widetilde{q}}(t,\theta)^{-1} .
\end{gather}

From here on, we will denote by $q$ the periodic potential that has been obtained from the complex curvature of $\gamma$ by the above regauging. We will omit the tilde in the names of
$\widetilde{q}$, $\widetilde{\lambda}$ and associated objects. Moreover we will henceforth usually omit the superscript ${}^q$ in $\alpha^q$ and similar objects.

Even where the potential $q$ is periodic, the extended frame $F$ is generally not periodic. The extent of its non-periodicity is measured by the \emph{monodromy}
$M(\lambda) := F(T,\lambda)$.

\begin{Remark} \label{th:monodromy-remarks}
The monodromy inherits the following properties from $F$.
\begin{enumerate}\itemsep=0pt
\item[{\rm{(i)}}] From the ODE \eqref{eq:frame1} it follows that the map $\bbC \to \mathrm{SL}_2(\bbC)$, $\lambda \mapsto M(\lambda)$ is analytic.

\item[{\rm{(ii)}}] From the reality condition \eqref{eq:reality} it follows that
\begin{gather} \label{eq:mreality}
	\overline{M\big(\bar\lambda\big)}^t = M(\lambda)^{-1} \qquad \mbox{for all} \quad \lambda \in \bbC.
\end{gather}
In particular $M(\lambda) \in \mathrm{SU}_2$ for all $\lambda \in \bbR$.
\item[{\rm{(iii)}}] The trace $\lambda \mapsto \operatorname{tr} M(\lambda)$ is analytic and satisfies
\begin{gather*} %\label{eq:trace}
	\operatorname{tr} \overline{M\big(\bar\lambda\big)} = \operatorname{tr} M(\lambda).
\end{gather*}
\item[{\rm{(iv)}}] The two eigenvalues $\bbC \to \bbC^\times,\lambda \mapsto \mu(\lambda)^{\pm 1}$ of $M(\lambda)$ are given by
\begin{gather*} %\label{eq:eigenvalues}
	\mu(\lambda)^{\pm 1} = \tfrac{1}{2} \Bigl( \operatorname{tr}M(\lambda) \pm \sqrt{\operatorname{tr}^2M(\lambda) - 4} \Bigr)
\end{gather*}
and satisfy $\operatorname{tr}M(\lambda) = \mu(\lambda) + \mu(\lambda)^{-1}$ and
\begin{gather*} %\label{eq:mu-reality}
	\overline{\mu\big(\bar{\lambda}\big)} = \mu(\lambda)^{-1}.
\end{gather*}
\item[{\rm{(v)}}] The map $\lambda \mapsto \mu(\lambda)$ is branched at the roots of $\operatorname{tr}^2M(\lambda) - 4$ of odd order, and there
\begin{gather*} %\label{eq:branch}
	\operatorname{tr}M(\lambda) = \pm 2 \Longleftrightarrow \mu(\lambda) + \mu(\lambda)^{-1} = \pm 2 \Longleftrightarrow \mu(\lambda) = \pm 1.
\end{gather*}
\item[{\rm{(vi)}}] The curve defined by equation~\eqref{eq:gamma-S3} or \eqref{eq:gamma-R3} is $T$-periodic if and only if
\begin{gather*} %\label{eq:closing1}
\bbS^3\colon \ M(1+\theta) = M(-1+\theta) = \pm \mathbbm{1} ; \qquad \bbR^3 \colon \ M(\theta) = \pm \mathbbm{1}, \qquad M'(\theta) = 0 .
\end{gather*}
\end{enumerate}
\end{Remark}

\begin{Lemma} \label{thm:series}
Let $\alpha\colon [0,T] \to \C^{2 \times 2}$ be Lebesgue integrable. Then the map
\begin{gather} \label{eq:series}
 t \mapsto \one + \sum_{n=1}^\infty \int_{t_n=0}^t \int_{t_{n-1}=0}^{t_n} \cdots
 \int_{t_1=0}^{t_2} \alpha(t_1)\alpha(t_2)\cdots \alpha(t_n)
 {\rm d}t_1\cdots {\rm d}t_n
\end{gather}
converges to the solution of $\tfrac{{\rm d}F}{{\rm d}t} = F\alpha$
with $F(0) = \one$.
\end{Lemma}

\begin{proof} The series in \eqref{eq:series} converges, because
\begin{gather*}
\left\| \int_0^t \int_0^{t_n} \cdots
 \int_0^{t_2} \alpha(t_1)\alpha(t_2)\cdots \alpha(t_n)
 {\rm d}t_1\cdots {\rm d}t_n \right\| \\
\qquad{} \leq \int_0^t \int_0^{t_n} \cdots
 \int_0^{t_2} \| \alpha(t_1)\|\|\alpha(t_2)\|\cdots
 \|\alpha(t_n)\|
 {\rm d}t_1\cdots {\rm d}t_n \\
\qquad{} \leq \frac{1}{n!}\int_0^t \int_0^t \cdots
 \int_0^t \|\alpha(t_1)\|\|\alpha(t_2)\|\cdots
 \|\alpha(t_n)\|
 {\rm d}t_1\cdots {\rm d}t_n \leq
\frac{1}{n!}\left(\int_0^t \| \alpha(s)\|{\rm d}s \right)^n.\tag*{\qed}
\end{gather*}\renewcommand{\qed}{}
\end{proof}

\section{Asymptotic analysis of the monodromy}
\label{Se:asymptotic}
The following theorem provides asymptotic estimates for the monodromy of a periodic poten\-tial~$q$.
We will see in Section~\ref{Se:perturbed-fourier} that this result permits us to derive the asymptotic behaviour of the perturbed Fourier coefficients of the potential of a given (periodic) curve,
and for this reason the theorem of the present section is fundamental for our treatment of curves.
The methods by which these asymptotic estimates are obtained were developed in \cite{klein-habil} for the sinh-Gordon equation. More specifically, Theorem~\ref{thm:asymp}(i) is analogous to the ``basic asymptotic'' of \cite[Theorem~5.4]{klein-habil}, and Theorem~\ref{thm:asymp}(ii) is analogous to the ``Fourier asymptotic'' of \cite[Theorem~7.1]{klein-habil}. Compared to the proofs in \cite{klein-habil}, the proof of Theorem~\ref{thm:asymp} is significantly simplified by the far simpler structure of our $\alpha$, see equation~\eqref{eq:alpha_q}, compared to the corresponding connection form for the sinh-Gordon equation. Note that the results of \cite{klein-habil} are summarized in \cite{MR3693920}.

In the sequel, we denote by $|\dots|$ also the maximum absolute row sum norm for $(2\times 2)$-matrices.

{\bf The vacuum.} The simplest curves in the space forms are the geodesics; they correspond to the \emph{vacuum} potential $q \equiv 0$ and up to isometries have an extended frame
\begin{gather}\label{thm:asymp:eq:F0}
	F_0 (t,\lambda) = \exp \big( \tfrac{1}{2} \lambda t \varepsilon \big)
\end{gather}
with monodromy
\begin{gather} \label{eq:M0}
	M_0(\lambda) = \exp \big( \tfrac{1}{2} \lambda T \varepsilon \big) .
\end{gather}
Now $M_0(\lambda) = \pm \mathbbm{1}$ if and only if
\begin{gather} \label{eq:lambdak0}
\lambda=\lambda_{k,0} := \tfrac{2\pi}{T}k \qquad \text{holds for some $k\in \bbZ$,}
\end{gather}
and then $M_0(\lambda_{k,0})=(-1)^k \one$. Also $|M_0(\lambda)| = \big| \exp\big( \tfrac{\mi T}{2}\lambda \big) \big| + \big| \exp\big( {-}\tfrac{\mi T}{2}\lambda \big) \big|$ and thus
\begin{gather} \label{eq:M0-estim}
	\exp\left( \tfrac{T}{2}|\IM(\lambda)| \right) \leq |M_0(\lambda)| \leq 2\exp\left( \tfrac{T}{2}|\IM(\lambda)| \right) .
\end{gather}
\begin{Theorem} \label{thm:asymp}
Let $q\in L^2([0,T])$ be a periodic potential with associated monodromy $M(\lambda)$. In comparison with the vacuum monodromy $M_0(\lambda)$ we then have the following:
\begin{enumerate}\itemsep=0pt
\item[{\rm (i)}] For every $\varepsilon>0$ there exists $R>0$ such that for any $\lambda \in \C$ with $|\lambda|\geq R$ we have
\begin{gather} \label{thm:asymp:eq:basic}
	|M(\lambda)-M_0(\lambda)| \leq \varepsilon |M_0(\lambda)| .
\end{gather}
\item[{\rm (ii)}] We have
\begin{gather*}
 |M(\lambda_{k,0})-M_0(\lambda_{k,0})| \in \ell^2(k) .
\end{gather*}
\end{enumerate}
Both estimates hold uniformly if $q$ varies over a relatively compact subset of $L^2([0,T])$.
\end{Theorem}

\begin{proof}Let $q\in L^2([0,T])$ be extended to $\bbR$ periodically, and $\alpha^q$ be as in~\eqref{eq:alpha_q}. Let $F$ solve~\eqref{eq:frame1}, and set $M(\lambda)=F(T,\lambda)$. The vaccuum extended frame is denoted as in~\eqref{thm:asymp:eq:F0} by $F_0(t,\lambda) = \exp (t \alpha_0)$ where $\alpha_0 = \tfrac{1}{2}\lambda\varepsilon$. Note that $F_0(\cdot,\lambda)\colon \R \to \mathrm{SL}(2,\C)$ is a homomorphism of Lie groups, so that for all $s,t \in \R$
\begin{gather}\label{thm:asymp:eq:F0-homo}
F_0(t,\lambda) F_0(s,\lambda) = F_0(t+s,\lambda) \qquad\text{and}\qquad F_0(t,\lambda)^{-1} = F_0(-t,\lambda).
\end{gather}
For $M_0(\lambda)=F_0(T,\lambda)$ and $t\in \R$ we have
\begin{gather}\label{thm:asymp:eq:F0-norm}
\exp\left( \tfrac12|\IM(\lambda)||t| \right) \leq |F_0(t,\lambda)| \leq 2 \exp\left( \tfrac12|\IM(\lambda)||t| \right).
\end{gather}
To prove the inequalities we need to estimate $(M(\lambda)-M_0(\lambda)) M_0(\lambda)^{-1}
= M(\lambda) M_0(\lambda)^{-1}-\one$. To do so, we develop $E(t,\lambda) := F(t,\lambda) F_0^{-1} (t,\lambda)$ as a power series with respect to~$q$. We will use this power series expansion to write $\lambda \mapsto E(T,\lambda)-\one$ explicitly as the Fourier transform of an $L^2$-function, see~\eqref{thm:asymp:eq:E-fourier} below. We are then able to derive the desired estimates by applying results from the theory of Fourier transforms in this situation.

We write $\alpha^q = \alpha_0 + \beta$ where $\beta = \tfrac{1}{2}( q\varepsilon_+ + \bar{q}\varepsilon_- )$. Note that $\beta$ is $\lambda$-independent and that $\alpha_0 \beta = -\beta \alpha_0$ holds; as a consequence of the latter equation, we have for every $t \in \R$, $s\in [0,T]$
\begin{gather}\label{thm:asymp:eq:beta-F0}
\beta(t) F_0(s,\lambda) = F_0(-s,\lambda) \beta(t) .
\end{gather}
Now $\tfrac{{\rm d}E}{{\rm d}t} = F\beta F_0^{-1} = E \widetilde{\beta}$ with $\widetilde{\beta} := F_0\beta F_0^{-1}$ and $ E(0,\lambda) = F(0,\lambda) F_0 (0,\lambda)^{-1} = \one$. By Lemma~\ref{thm:series} (with $\widetilde{\beta}$ in the place of $\alpha$) we obtain
\begin{gather} \label{thm:asymp:eq:E}
	M(\lambda)M_0(\lambda)^{-1} = E(T,\lambda) = \one + \sum_{n=1}^\infty E_n(\lambda)
\end{gather}
with
\begin{gather}
E_n(\lambda) = \int_{t_n=0}^T \int_{t_{n-1}=0}^{t_n} \cdots \int_{t_1=0}^{t_2} \widetilde{\beta}(t_1,\lambda)\widetilde{\beta}(t_2,\lambda)\cdots \widetilde{\beta}(t_n,\lambda)
\mathrm{d}^nt \nonumber\\
\hphantom{E_n(\lambda)}{} = \int_{t_n=0}^T \int_{t_{n-1}=0}^{t_n} \cdots \int_{t_1=0}^{t_2} \beta(t_1)\beta(t_2)\cdots \beta(t_n)F_0(\xi_n(t),\lambda)
\mathrm{d}^nt\label{thm:asymp:eq:En}
\end{gather}
and
$\xi_n(t) = 2(-t_n+t_{n-1}-t_{n-2}-\dots+(-1)^{n}t_1)$, where the second equals sign in~\eqref{thm:asymp:eq:En} follows from equations~\eqref{thm:asymp:eq:F0-homo} and~\eqref{thm:asymp:eq:beta-F0}.

In the integral of equation~\eqref{thm:asymp:eq:En} we now carry out the substitution $(t_1,\dots,t_n) \mapsto (s_1,\dots,s_n)$ with $s_1=t_1$ and $s_j = t_j-s_{j-1}$ for $2 \leq j\leq n$, so that
\begin{gather*}
s_1 = t_1,\qquad
s_2 = t_2-s_1 = t_2-t_1,\qquad s_3 = t_3-s_2 = t_3-t_2+t_1, \qquad
\dots, \\
s_n = t_n-s_{n-1} = t_n-t_{n-1}+t_{n-2}- \dots + (-1)^{n+1}t_1 = -\tfrac12\xi_n(t) .
\end{gather*}
Then we have $t_1=s_1$ and $t_j = s_j+s_{j-1}$ for $2 \leq j\leq n$.
Thus, the corresponding mapping $\Phi\colon (s_j) \mapsto (t_j)$ is a diffeomorphism with $\det \Phi'=1$ from
\begin{gather*}
 U := \big\{(s_1,\dots,s_n)\in \R^n\,|\,\forall\, j=1,\dots,n\colon s_j\geq 0,\, s_{n-1}+s_n \leq T,\, \forall\, j=3,\dots,n\colon s_{j-2} \leq s_{j}\big\}
\end{gather*}
onto the simplex $ \big\{(t_1,\dots,t_n)\in \R^n\,|\,0 \leq t_1 \leq t_{2} \leq \cdots \leq t_n \leq T\big\}$ which is the domain of integration in~\eqref{thm:asymp:eq:En}. This substitution into~\eqref{thm:asymp:eq:En} yields
\begin{gather}
E_n(\lambda) = \int_{s_n=0}^T \int_{s_{n-1}=0}^{T-s_n} \int_{s_{n-2}=0}^{s_n} \int_{s_{n-3}=0}^{s_{n-1}} \cdots \nonumber\\
\hphantom{E_n(\lambda) =}{} \cdots \int_{s_2=0}^{s_4} \int_{s_1=0}^{s_3}
\beta(s_1) \beta(s_1+s_2) \beta(s_2+s_3) \cdots \beta(s_{n-1}+s_n)F_0(-2s_n,\lambda)\mathrm{d}^ns\nonumber \\
\hphantom{E_n(\lambda)}{} = \int_0^T G_n(s)F_0(-2s,\lambda)\mathrm{d}s\label{thm:asymp:eq:EnGn}
\end{gather}
with $G_1 := \beta$ and for $n\geq 2$
\begin{gather*}
G_n(s_n) = \int_{s_{n-1}=0}^{T-s_n} \int_{s_{n-2}=0}^{s_n} \int_{s_{n-3}=0}^{s_{n-1}} \cdots \\ \hphantom{G_n(s_n) =}{} \cdots \int_{s_2=0}^{s_4} \int_{s_1=0}^{s_3}
\beta(s_1) \beta(s_1+s_2) \beta(s_2+s_3) \cdots \beta(s_{n-2}+s_{n-1})\beta(s_{n-1}+s_n)\mathrm{d}^{n-1}s .
\end{gather*}
Note that $G_n(s_n)$ is well-defined for every $s_n \in [0,T]$. We will now show that $G_n$ is bounded for each $n\geq 2$ and in fact obtain an explicit upper bound for $|G_n(s)|$. For even $n\geq 2$, we fix $s_n\in [0,T]$ and define the simplices
\begin{gather*}
S_+ := \big\{(s_2,s_4,s_6,\dots,s_{n-2})\in\R^{n/2-1}\,|\,0 \leq s_2 \leq s_4 \leq \dots \leq s_{n-2} \leq s_n \big\} , \\
S_- := \big\{(s_1,s_3,s_5,\dots,s_{n-1})\in \R^{n/2}\,|\,0 \leq s_{1} \leq s_3 \dots \leq s_{n-1} \leq T-s_n\big\} .
\end{gather*}
Then we have
\begin{gather*}
 G_n(s_n) = \int_{S_+} \int_{S_-} \beta(s_1) \beta(s_1+s_2) \beta(s_2+s_3) \cdots \beta(s_{n-1}+s_n) \mathrm{d}^{n/2}s\mathrm{d}^{n/2-1}s .
\end{gather*}
For the inner integral, we have by the Cauchy--Schwarz inequality
\begin{gather*}
 \left | \int_{S_-} \beta(s_1) \beta(s_1+s_2) \beta(s_2+s_3) \cdots \beta(s_{n-1}+s_n)\mathrm{d}^{n/2}s \right| \leq \|\beta\|_{L^2}^n ,
\end{gather*}
and thus we have
\begin{gather*}
 |G_n(s_n)| \leq \mathrm{vol}(S_+) \|\beta\|_{L^2}^n = \frac{s_n^{n/2-1}}{(n/2-1)!} \|\beta\|_{L^2}^n \leq \frac{T^{n/2-1}}{(n/2-1)!} \|\beta\|_{L^2}^n .
\end{gather*}
For odd $n\geq 3$ we argue similarly and conclude that
\begin{gather*}
 |G_n(s_n)| \leq \frac{T^{(n-1)/2}}{((n-1)/2)!} \|\beta\|_{L^2}^n .
\end{gather*}
Thus we obtain for integers $n\geq 2$ of any parity
\begin{gather} \label{thm:asymp:eq:Gn-estimate}
 |G_n(s_n)| \leq \frac{T^{\lfloor n/2 \rfloor-1}}{(\lfloor n/2 \rfloor-1)!} \|\beta\|_{L^2}^n .
\end{gather}
Here $\lfloor x \rfloor$ denotes the greatest integer that is less than or equal to $x$.
The estimate~\eqref{thm:asymp:eq:Gn-estimate} shows that $G_n$ is bounded on $[0,T]$ for every $n\geq 2$, and therefore in particular $G_n \in L^2\big([0,T],\C^{2\times 2}\big)$. We also have $G_1 = \beta \in L^2\big([0,T],\C^{2\times 2}\big)$. The estimate~\eqref{thm:asymp:eq:Gn-estimate} further shows that $\sum\limits_{n=2}^\infty G_n$ converges uniformly on $[0,T]$ (because the power series $r \mapsto \sum\limits_{n=1}^\infty \frac{T^{\lfloor n/2 \rfloor-1}}{(\lfloor n/2 \rfloor-1)!} r^n$ is convergent for every real~$r$) and therefore in particular in $L^2\big([0,T],\C^{2\times 2}\big)$, and thus
\begin{gather*}
 G := \sum_{n=1}^\infty G_n
\end{gather*}
converges to a function $G \in L^2\big([0,T],\C^{2\times 2}\big)$. It follows from equations~\eqref{thm:asymp:eq:E} and \eqref{thm:asymp:eq:EnGn} that
\begin{gather} \label{thm:asymp:eq:E-fourier}
	M(\lambda)M_0(\lambda)^{-1} -\one = \int_0^T G(s)F_0(-2s,\lambda)\mathrm{d}s,
\end{gather}
and by equation~\eqref{thm:asymp:eq:F0}, the entries of this matrix are of the form
\begin{gather} \label{thm:asymp:eq:E-form}
	\int_0^T g_{\mu\nu}(s)\exp(\mp \mi \lambda s)\mathrm{d}s
\end{gather}
with functions $g_{\mu\nu}\in L^2([0,T])$. This shows that for real $\lambda$, the map $\lambda \mapsto M(\lambda) M_0(\lambda)^{-1}-\one$ is essentially the Fourier transform of the $L^2$-function~$G(s)$.

We first prove part (ii) of the theorem. For $\lambda=\lambda_{k,0}=\tfrac{2\pi k}{T}$ we have $M_0(\lambda_{k,0}) = (-1)^k \one$ and therefore by equation~\eqref{thm:asymp:eq:E-fourier}
\begin{gather*}
 M(\lambda_{k,0}) - M_0(\lambda_{k,0}) = \bigr( M(\lambda) M_0(\lambda)^{-1} -\one \bigr) M_0(\lambda_{k,0}) = (-1)^k\int_0^T G(s)F_0(-2s,\lambda)\mathrm{d}s .
\end{gather*}
The entries of the right-hand matrix are of the form
\begin{gather*}
 (-1)^k \int_0^T g_{\mu\nu}(s)\exp(\mp \mi \lambda_{k,0} s)\mathrm{d}s = (-1)^k \int_0^T g_{\mu\nu}(s)\exp\big(\mp \tfrac{2\pi\mi k}{T} s\big)\mathrm{d}s
\end{gather*}
with functions $g_{\mu\nu}\in L^2([0,T])$, and therefore sequences in $\ell^2(k)$ by Plancherel's theorem (see, e.g., \cite[Proposition~3.2.7(1)]{grafakos}). This implies part (ii) of the theorem.

For real values of $\lambda$, part (i) of the theorem follows similarly via the Riemann--Lebesgue lemma (see, e.g., \cite[Proposition~2.2.17]{grafakos}).
But because $\lambda$ can be complex-valued, we need to apply a more refined argument. Let $\varepsilon>0$ be given. We choose $\delta > 0$ at first arbitrarily.
Then by equations~\eqref{thm:asymp:eq:E-fourier} and~\eqref{thm:asymp:eq:F0-homo} we have
\begin{gather*}
M(\lambda)-M_0(\lambda) = (M(\lambda)M_0(\lambda)^{-1}-\one) M_0(\lambda) \\
\hphantom{M(\lambda)-M_0(\lambda)}{} = \int_0^T G(s)F_0(-2s,\lambda)\mathrm{d}s F_0(T,\lambda) = \int_0^T G(s)F_0(T-2s,\lambda)\mathrm{d}s \\
\hphantom{M(\lambda)-M_0(\lambda)}{} = \int_0^\delta G(s)F_0(T-2s,\lambda)\mathrm{d}s + \int_\delta^{T-\delta} G(s)F_0(T-2s,\lambda)\mathrm{d}s \\
\hphantom{M(\lambda)-M_0(\lambda)=}{} + \int_{T-\delta}^T G(s)F_0(T-2s,\lambda)\mathrm{d}s .
\end{gather*}
For $s \in [0,\delta] \cup [T-\delta,T]$ we have $|F_0(T-2s,\lambda)| \leq |F_0(T,\lambda)| = |M_0(\lambda)|$, and for $s \in [\delta,T-\delta]$ we have by~\eqref{thm:asymp:eq:F0-norm}
\begin{gather*}
 |F_0(T-2s,\lambda)| \leq | F_0(T-2\delta,\lambda)| \leq 2\left| \exp\left( \tfrac12 |\IM(\lambda)|(T-2\delta) \right)\right| \leq 2{\rm e}^{-|\IM(\lambda)|\delta} |M_0(\lambda)| .
\end{gather*}
Therefore we obtain that
\begin{gather*}
| M(\lambda)-M_0(\lambda) | \leq \int_0^\delta |G(s)||F_0(T-2s,\lambda)|\mathrm{d}s + \int_\delta^{T-\delta} |G(s)| |F_0(T-2s,\lambda)|\mathrm{d}s \\
\hphantom{| M(\lambda)-M_0(\lambda) | \leq}{} + \int_{T-\delta}^T |G(s)||F_0(T-2s,\lambda)|\mathrm{d}s \\
\hphantom{| M(\lambda)-M_0(\lambda) |}{} \leq \| G|[0,\delta] \|_{L^1} |M_0(\lambda)| + \| G \|_{L^1} \cdot 2{\rm e}^{-|\IM(\lambda)|\delta} |M_0(\lambda)| \\
\hphantom{| M(\lambda)-M_0(\lambda) | \leq}{}+ \| G|[T-\delta,T] \|_{L^1} |M_0(\lambda)| \\
\hphantom{| M(\lambda)-M_0(\lambda) |}{} = \big( \| G|[0,\delta] \|_{L^1} + \| G \|_{L^1} \cdot 2{\rm e}^{-|\IM(\lambda)|\delta} + \| G|[T-\delta,T] \|_{L^1} \big) |M_0(\lambda)| .
\end{gather*}
We now choose $\delta>0$ so that $ \| G|[0,\delta]\|_{L^1}, \| G|[T-\delta,T] \|_{L^1} \leq \tfrac{\varepsilon}{4}$ holds. Then we choose $C>0$ so that $ \| G \|_{L^1}\cdot 2{\rm e}^{-C \delta} \leq \tfrac{\varepsilon}{2}$ holds. With these choices, \eqref{thm:asymp:eq:basic} holds for all $\lambda\in\C$
with $|\IM(\lambda)|\geq C$.

It remains to show that~\eqref{thm:asymp:eq:basic} also holds within the horizontal strip $\bigr\{\lambda\in\C \,|\, |\IM(\lambda)|\leq C\bigr\}$ for~$\lambda$ of sufficiently
large modulus. Because of
\begin{gather*}
 | M(\lambda)-M_0(\lambda) | \leq \big|M(\lambda)M_0(\lambda)^{-1}-\one| |M_0(\lambda)\big|
\end{gather*}
and the entries of $M(\lambda)M_0(\lambda)^{-1}-\one$ are of the form given in~\eqref{thm:asymp:eq:E-form}, this statement immediately follows from the variant of the Riemann--Lebesgue
lemma given as Lemma~\ref{lem:riemannlebesgue} below.
\end{proof}

The following variant of the Riemann--Lebesgue lemma was used in the preceding proof.
\begin{Lemma}\label{lem:riemannlebesgue} Let $N$ be a compact subset of $L^1([a,b])$. For any given $\varepsilon,C>0$ there then exists $R>0$ such that for every $\lambda\in \C$ with $|\IM(\lambda)|\leq C$, $|\RE(\lambda)|\geq R$ and every $g\in N$ we have
\begin{gather*}
 \left| \int_a^b g(t){\rm e}^{-2\pi \mi\lambda t}\mathrm{d}t \right| \leq \varepsilon .
\end{gather*}
\end{Lemma}

\begin{proof} We extend the functions in $L^1([a,b])$ to $\R$ by zero, then we have $N \subset L^1(\R)$. The map
\begin{gather*}
 \Phi\colon \ L^1(\R)\times [-C,C] \to L^1(\R), \qquad (g,y) \mapsto \big(t \mapsto {\rm e}^{2\pi yt} g(t)\big)
\end{gather*}
is continuous, hence the image $\widetilde{N} := \Phi(N\times [-C,C])$ is a compact set in $L^1(\R)$.

Therefore there exist finitely many $f_1,\dots,f_n\in L^1(\R)$ so that $\widetilde{N} \subset \bigcup_{k=1}^n U^{L^1}_{\varepsilon/2}(f_k)$. By the classical Riemann--Lebesgue lemma \cite[Proposition~2.2.17]{grafakos}, there exists $R>0$, so that we have $\big|\widehat{f}_k(x)\big|\leq\tfrac{\varepsilon}{2}$ for every $x\in \R$ with $|x|\geq R$ and every $k\in \{1,\dots,n\}$; here we denote for any $f\in L^1(\R)$ by
\begin{gather*}
 \widehat{f}(x) := \int_{\R} f(t){\rm e}^{-2\pi \mi xt}\mathrm{d}t
\end{gather*}
the Fourier transform of $f$.

Now let $g\in N$ and $\lambda=x+\mi y \in \C$ be given with $|x| \geq R$, $|y| \leq C$. By construction,
there exists some $k\in \{1,\dots,n\}$ with $\|\Phi(g,y)-f_k\|_1 \leq \tfrac{\varepsilon}{2}$,
and with this $k$ we have
\begin{gather*}
 \big\|\widehat{\Phi(g,y)} - \widehat{f}_k\big\|_\infty \leq \|\Phi(g,y)-f_k\|_1 \leq \frac{\varepsilon}{2}
\end{gather*}
by the Hausdorff--Young inequality for the case $p'=\infty$, $p=1$ \cite[Proposition~2.2.16]{grafakos}. We now have
\begin{gather*}
\left| \int_a^b g(t){\rm e}^{-2\pi \mi \lambda t}\mathrm{d}t \right|
= \left| \int_a^b g(t){\rm e}^{-2\pi \mi x t}{\rm e}^{2\pi y t}\mathrm{d}t \right|
= \left| \int_a^b \Phi(g,y)(t){\rm e}^{-2\pi \mi x t}\mathrm{d}t \right| = \big| \widehat{\Phi(g,y)}(x) \big| \\
\hphantom{\left| \int_a^b g(t){\rm e}^{-2\pi \mi \lambda t}\mathrm{d}t \right|}{}
\leq \big|\widehat{\Phi(g,y)}(x) - \widehat{f}_k(x)\big| + \big| \widehat{f}_k(x) \big|
\leq \frac{\varepsilon}{2} + \frac{\varepsilon}{2} = \varepsilon .\tag*{\qed}
\end{gather*}\renewcommand{\qed}{}
\end{proof}

\begin{Corollary}\label{cor:asymp-M'}Let $q\in L^2([0,T])$ be a periodic potential, and $M(\lambda)$ be the monodromy associated to $q$. For every $\varepsilon>0$ there exists $R>0$ so that for any $\lambda \in \C$ with $|\lambda|\geq R$ we have
\begin{gather*}
 |M'(\lambda)-M_0'(\lambda)| \leq \varepsilon |M_0(\lambda)| ,
\end{gather*}
where $'$ denotes the derivative with respect to $\lambda$. This estimate holds uniformly if $q$ varies over a relatively compact subset of $L^2([0,T])$.
\end{Corollary}

\begin{proof}We first note that because of \eqref{eq:M0-estim} there exists $C>0$ so that we have
\begin{gather*}
 \max_{\lambda'\in U_1(\lambda)} |M_0(\lambda')| \leq C |M_0(\lambda)| \qquad \text{for every $\lambda \in \bbC$.}
\end{gather*}
Now let $\varepsilon > 0$ be given. By Theorem~\ref{thm:asymp}(1) there exists $R_1>0$ so that we have
\begin{gather*}
 |M(\lambda)-M_0(\lambda)| \leq \frac{\varepsilon}{C} |M_0(\lambda)| \qquad \text{for every $\lambda\in \bbC$ with $|\lambda|\geq R_1$.}
\end{gather*}
Put $R := R_1+1$. For $\lambda\in \C$ with $|\lambda|\geq R$ we then have $\overline{U_1(\lambda)} \subset \{|\lambda|\geq R_1\}$ and therefore by Cauchy's Inequality, applied to the holomorphic function
$M-M_0$
\begin{gather*}
 |M'(\lambda)-M_0'(\lambda)| \leq \max_{\lambda' \in U_1(\lambda)} |M(\lambda')-M_0(\lambda')| \leq \frac{\varepsilon}{C} \max_{\lambda' \in U_1(\lambda)} |M_0(\lambda')| \leq \varepsilon |M_0(\lambda)| .\tag*{\qed}
\end{gather*}\renewcommand{\qed}{}
\end{proof}

\begin{Corollary}\label{cor:asymp} Let $q\in L^2([0,T])$ be a periodic potential, and $M(\lambda)$ be the monodromy associated to $q$. For every sequence $(\lambda_k)_{k\in \bbZ}$ with $\lambda_k-\lambda_{k,0} \in \ell^2(k)$ we then have
\begin{gather}\label{eq:cor:asymp:claim}
|M(\lambda_k)-M_0(\lambda_k)| \in \ell^2(k) .
\end{gather}
This estimate holds uniformly if $q$ varies over a relatively compact subset of $L^2([0,T])$.
\end{Corollary}

\begin{proof}We first note that there exists a horizontal strip $S$ around the $x$-axis in $\bbC$ so that the line $[\lambda_{k,0},\lambda_k]$ is contained in $S$ for all $k \in \bbZ$. $|M_0|$ is bounded on $S$, and hence it follows from Corollary~\ref{cor:asymp-M'} (applied with $\varepsilon=1$) that there exists $C,R>0$ so that we have
\begin{gather}\label{eq:cor:asymp:M'-asymp}
|M'(\lambda)-M_0'(\lambda)| \leq C \qquad \text{for $\lambda\in S$ with $|\lambda|\geq R$.}
\end{gather}
We now have for $k\in \bbZ$
\begin{gather}\label{eq:cor:asymp:MM0-equation}
M(\lambda_k)-M_0(\lambda_k) = M(\lambda_{k,0})-M_0(\lambda_{k,0}) + \int_{\lambda_{k,0}}^{\lambda_k} \left( M'(\lambda)-M_0'(\lambda)\right)\mathrm{d}\lambda .
\end{gather}
We have $M(\lambda_{k,0})-M_0(\lambda_{k,0}) \in \ell^2(k)$ by Theorem~\ref{thm:asymp}(2). Moreover we have for~$|k|$ sufficiently large by \eqref{eq:cor:asymp:M'-asymp}
\begin{gather*}
 \left| \int_{\lambda_{k,0}}^{\lambda_k} \left( M'(\lambda)-M_0'(\lambda)\right)\mathrm{d}\lambda \right| \leq C |\lambda_k-\lambda_{k,0}| ,
\end{gather*}
where the integration is carried out along the straight line from $\lambda_{k,0}$ to $\lambda_k$.
Because of $\lambda_k-\lambda_{k,0}\in \ell^2(k)$ it follows from this estimate that $\int_{\lambda_{k,0}}^{\lambda_k} \left( M'(\lambda)-M_0'(\lambda)\right)\mathrm{d}\lambda \in \ell^2(k)$ holds.
Thus \eqref{eq:cor:asymp:claim} follows from equation~\eqref{eq:cor:asymp:MM0-equation}.
\end{proof}

\section{Perturbed Fourier coefficients}\label{Se:perturbed-fourier}

Let a periodic potential $q\in L^2([0,T])$ be given, and let $M(\lambda)$ be the monodromy associated to it. We write
\begin{gather*}
 M(\lambda) = \begin{pmatrix} a(\lambda) & b(\lambda) \\ c(\lambda) & d(\lambda) \end{pmatrix}
\end{gather*}
with holomorphic functions $a,b,c,d\colon \C \to \C$.

\begin{Lemma}\label{lem:lambdak}The set of zeros (with multiplicities) of the function $a-d$ is enumerated by a~sequence $(\lambda_k)_{k\in \bbZ}$, such that $\lambda_k-\lambda_{k,0} \in \ell^2(k)$ holds.
\end{Lemma}

\begin{proof}We write the monodromy of the vacuum as $M_0(\lambda)=\left( \begin{smallmatrix} a_0(\lambda) & 0 \\ 0 & d_0(\lambda) \end{smallmatrix} \right)$ with $a_0(\lambda) = \exp\big( \tfrac{\mi T}{2}\lambda\big)$ and $d_0(\lambda) = \exp\big( {-}\tfrac{\mi T}{2}\lambda\big)$. The idea of the proof is to use Rouch\'e's theorem to compare the number of zeros of the function $f := a-d$ on suitable domains to the number of zeros of the function $f_0(\lambda) := a_0(\lambda)-d_0(\lambda)=2\sinh\big(\tfrac{\mi T}{2}\lambda\big)=2\mi\sin\big(\tfrac{T}{2}\lambda\big)$. Note that the zeros of the latter function are exactly the $\lambda_{k,0}$, $k\in \bbZ$.

We first note that there exist $\delta,C>0$ so that for any $\lambda\in \C$ with $|\lambda-\lambda_{k,0}|\geq \delta$ for all $k\in \bbZ$ we have
\begin{gather}\label{lem:lambdak:eq:sinh}
|f_0(\lambda)| = 2 \left| \sinh\left( \tfrac{\mi T}{2}\lambda \right)\right| \geq 2C \exp\left( \tfrac{T}{2}|\IM(\lambda)| \right) \overset{\eqref{eq:M0-estim}}{\geq} C |M_0(\lambda)| .
\end{gather}
Now let $\varepsilon := \tfrac12 C>0$. By Theorem~\ref{thm:asymp}(1) there exists $R>0$ so that for $\lambda\in \C$ with $|\lambda|\geq R$ we have
\begin{gather}\label{lem:lambdak:eq:asymp}
|f(\lambda)-f_0(\lambda)| \leq |M(\lambda)-M_0(\lambda)| \leq \varepsilon |M_0(\lambda)| .
\end{gather}
Then choose $K\in \bbN$ so that $\overline{U_\delta(\lambda_{k,0})} \subset \{|\lambda|\geq R\}$ holds for all $k\in \bbZ$ with $|k|\geq K$. For such $k$ and $\lambda \in \partial U_\delta(\lambda_{k,0})$
both estimates~\eqref{lem:lambdak:eq:sinh} and~\eqref{lem:lambdak:eq:asymp} apply, and thus we obtain
\begin{gather*}
 |f(\lambda)-f_0(\lambda)| \leq \varepsilon |M_0(\lambda)| \leq \frac{\varepsilon}{C} |f_0(\lambda)| = \frac12 |f_0(\lambda)| < |f_0(\lambda)| .
\end{gather*}
By Rouch\'e's theorem, it follows that the functions $f$ and $f_0$ have the same number of zeros on~$U_\delta(\lambda_{k,0})$. Because $\lambda_{k,0}$ is the only zero of~$f_0$ in~$U_\delta(\lambda_{k,0})$, $f$ also has exactly one zero in~$U_\delta(\lambda_{k,0})$, which we denote by~$\lambda_k$.

For $k\geq K$, we also consider the disk $V_k := U_{(\lambda_{k,0}+\lambda_{k+1,0})/2}(0)$. Then on $\partial V_k$, again both inequalities~\eqref{lem:lambdak:eq:sinh} and~\eqref{lem:lambdak:eq:asymp} apply, whence $|f(\lambda)-f_0(\lambda)|<|f_0(\lambda)|$ follows for $\lambda \in \partial V_k$. Again using Rouch\'e's theorem, we see that~$f$ and~$f_0$ have the same number of zeros on~$V_k$, namely $2k+1$ many. This argument shows that besides the points of $(\lambda_k)_{|k|> K}$, the function~$f$ has exactly $2K+1$ many zeros (counted with multiplicity), which we enumerate arbitrarily as $(\lambda_k)_{-K \leq k \leq K}$.

It remains to show that $\lambda_k-\lambda_{k,0}\in \ell^2(k)$ holds. For this we note that there exists $C>0$ so that we have for all $k\in \bbZ$
\begin{gather*}
 |f_0(\lambda_k)-f_0(\lambda_{k,0})| \geq C |\lambda_k-\lambda_{k,0}| .
\end{gather*}
Because of $f_0(\lambda_{k,0})=0=f(\lambda_k)$ we have
\begin{align*}
f_0(\lambda_k)-f_0(\lambda_{k,0}) & = f_0(\lambda_k)-f(\lambda_k) = f_0(\lambda_{k,0})-f(\lambda_{k,0}) - \int_{\lambda_{k,0}}^{\lambda_k} (f'(\lambda)-f_0'(\lambda))\mathrm{d}\lambda
\end{align*}
and therefore
\begin{gather*}
 |\lambda_k-\lambda_{k,0}| \leq \tfrac{1}{C} |f_0(\lambda_k)-f_0(\lambda_{k,0})| \\
 \hphantom{|\lambda_k-\lambda_{k,0}|}{}
 \leq \tfrac{1}{C} \Big( |f_0(\lambda_{k,0})-f(\lambda_{k,0})| + \max_{\lambda\in [\lambda_{k,0},\lambda_k]} |f'(\lambda)-f_0'(\lambda)| |\lambda_k-\lambda_{k,0}| \Big) .
\end{gather*}
We have $f_0(\lambda_{k,0})-f(\lambda_{k,0}) \in \ell^2(k)$ by Theorem~\ref{thm:asymp}(2), and moreover $\max\limits_{\lambda\in [\lambda_{k,0},\lambda_k]} |f'(\lambda)-f_0'(\lambda)| \leq \tfrac{C}{2}$ for all $k\in \bbZ$ with $|k|$ sufficiently large by Corollary~\ref{cor:asymp-M'}. Thus we obtain
\begin{gather*}
 |\lambda_k-\lambda_{k,0}| \leq \ell^2(k) + \tfrac12 |\lambda_k-\lambda_{k,0}|
\end{gather*}
and therefore $\lambda_k-\lambda_{k,0} \in \ell^2(k)$.
\end{proof}

\begin{Definition}\label{def:perturbedfourier} Let $(\lambda_k)_{k\in\bbZ}$ be as in Lemma~\ref{lem:lambdak}. Then we call the sequence $(z_k)_{k\in \bbZ}$ with $z_k := 2(-1)^kb(\lambda_k)$ the \emph{perturbed Fourier coefficients} of the potential $q$.
\end{Definition}

The proof of Theorem~\ref{thm:asymp} shows that the perturbed Fourier coefficients $(z_k)$ are asymptotically close to the usual Fourier coefficients of the potential~$q$, also see Lemma~\ref{lem:delta-zk} below. This is the reason for the name ``perturbed Fourier coefficients''.

\begin{Proposition}\label{prop:perturbedfourier}We have $z_k \in \ell^2(k)$.
\end{Proposition}

\begin{proof}We have $z_k = 2(-1)^kb(\lambda_k) = 2(-1)^k(b(\lambda_k)-b_0(\lambda_k))$ because of $b_0=0$, and therefore it follows from Corollary~\ref{cor:asymp} and Lemma~\ref{lem:lambdak} that $z_k \in \ell^2(k)$ holds.
\end{proof}

{\bf The spectral curve.} The spectral curve of $M(\lambda)$ or of the corresponding periodic poten\-tial~$q$ is the hyperelliptic complex curve defined by the characteristic equation of $M(\lambda)$
\begin{gather}\label{prop:finitegap-periodic:spectralcurve}
\Sigma = \bigr\{ (\lambda,\mu)\in \C^2 \,|\, \mu^2 - \Delta(\lambda)\mu+1=0 \bigr\} \qquad\text{with}\quad \Delta(\lambda) := \operatorname{tr}(M(\lambda)) ,
\end{gather}
compare for example \cite[Section~2.2, Example~3]{grinevich-schmidt1995} or \cite[Section~3,~1)]{Grinevich/Santini}. A holomorphic involution of $\Sigma$ is given by $\sigma\colon (\lambda,\mu) \mapsto \big(\lambda,\mu^{-1}\big)$, in this sense $\Sigma$ is hyperelliptic above~$\bbC$,
branch points of~$\Sigma$ occur at the zeros of~$\Delta^2-4$ of odd order, singularities of $\Sigma$ occur at the zeros of $\Delta^2-4$ of order $\geq 2$, and the eigenbundle of $M(\lambda)$ is a holomorphic line bundle $\Lambda$ on a suitable partial desingularisation $\widehat{\Sigma}$ of $\Sigma$ (see \cite[Section~4]{klss2016}, where $\widehat{\Sigma}$ is called the $\mathcal{S}$-halfway normalisation of the holomorphic matrix $M(\lambda)$). In general, $\widehat{\Sigma}$ can have infinite geometric genus, however if we regard $\widehat{\Sigma}$ as a complex space, this complex space can be compactified by adding points above $\lambda=\infty$ with a special topology described in \cite[Chapter~2]{schmidt1996}.
The compactified surface, which we again denote by $\widehat{\Sigma}$, is a hyperelliptic surface above $\CPone$ with the hyperelliptic involution~$\sigma$. For~$\lambda\to\infty$, $\widehat{\Sigma}$~is approximated by the spectral curve of the vacuum $q=0$, which shows that~$\infty$ is not a branching point of the compactification $\widehat{\Sigma}$, in other words there are two points $\infty_+,\infty_- \in \widehat{\Sigma}$ that are above $\infty \in \CPone$.

Recall that $q$ is said to have \emph{finite gaps} if $q$ is stationary under all but finitely many flows of the NLS hierarchy. This is the case if and only if $\widehat{\Sigma}$
has finite arithmetic genus, see \cite[Section~2]{calini-ivey2005}. If this is the case, then the compactification $\widehat{\Sigma}$ described above is the compactification
in the usual sense.

\begin{Proposition}
\label{prop:fourier-finitegap}
The potential $q \in L^2([0,T])$ is finite gap if and only if all but finitely many of the perturbed Fourier coefficients $z_k$ vanish.
\end{Proposition}

\begin{proof}
We first show that in any event, $\lambda_k\in \bbR$ holds for all but finitely many $k\in \bbZ$.
If some $\lambda\in\bbC$ is a zero of $a-d$, then $\overline{\lambda}$ also is a zero of $a-d$ because of equation~\eqref{eq:mreality}. On the other hand,
we know from Lemma~\ref{lem:lambdak} that $U_\delta(\lambda_{k,0})$ (where $\delta>0$ is small) contains exactly one zero of $a-d$ for all but finitely many $k\in\bbZ$.
Due to $\lambda_{k,0}\in\bbR$, the relation $\lambda_k \in U_\delta(\lambda_{k,0})$ implies $\overline{\lambda_k}\in U_\delta(\lambda_{k,0})$, and therefore $\overline{\lambda_k} = \lambda_k$
holds for all but finitely many $k\in\bbZ$.

Whenever $\lambda_k\in\bbR$ holds, we have $M(\lambda_k) \in \mathrm{SU}_2$ by Remark~\ref{th:monodromy-remarks}(ii), and therefore $0=\tfrac{(-1)^k}{2}z_k=b(\lambda_k)$ implies
$M(\lambda_k)=\pm\one$. Therefore $\lambda_k$ then is a zero of $\Delta^2-4$ of order $\geq 2$, hence a double point of the spectral curve $\Sigma$ of $q$, and moreover $\lambda_k$
is in the support of the spectral divisor of $q$.
In this case $\widehat{\Sigma}$ is at $\lambda=\lambda_k$ the normalisation of $\Sigma$, hence $\widehat{\Sigma}$ has neither singularities nor branch points here.
This shows that if all but finitely many of the $z_k$ vanish, then $\widehat{\Sigma}$ has finite arithmetic genus.

Conversely, if $\widehat{\Sigma}$ has finite artihmetic genus, then all but finitely many of the zeros of $\Delta^2-4$ are double points of $\Sigma$ that are in the
support of the spectral divisor of $q$. If this is the case for some $\lambda$, then $M(\lambda) = \pm \one$ holds,
and therefore $\lambda=\lambda_k$ for some $k\in \bbZ$ and $z_k=0$.
\end{proof}

%
%%%%%%%%%%%%%%%%%%%%%%%%%%%%%%
%%%%%%%%%%%%%%%%%%%%%%%%%%%%%%
%

\section{Asymptotic analysis of the variation of the monodromy}
\label{Se:var-monodromy}

In Section~\ref{Se:perturbed-fourier-map}
we will construct finite gap potentials by means of certain perturbed Fourier coefficient maps on the space of potentials. To show that such maps are invertible, we need to consider
variations of the Fourier coefficents corresponding to a variation of the potential $q$. To facilitate this investigation, we begin in the present section by studying variations of the monodromy.

\begin{Lemma}
\label{lem:asymp-F}
Let $q\in L^2([0,T])$ be a periodic potential with the extended frame $F=F(t,\lambda)$. Let $F_0=F_0(t,\lambda)$ be the extended frame for the vacuum potential $q=0$.
\begin{enumerate}\itemsep=0pt
\item[{\rm{(i)}}] For every $\varepsilon>0$ there exists $R>0$ so that for any $\lambda\in\bbC$ with $|\lambda|\geq R$ and any $t\in [0,T]$ we have
\begin{gather*}%\label{lem:asymp-F1}
	|F(t,\lambda)-F_0(t,\lambda)| \leq \varepsilon |F_0(t,\lambda)| \leq \varepsilon |M_0(\lambda)| . \end{gather*}
\item[{\rm{(ii)}}] Let a sequence $(\lambda_k)_{k\in \bbZ}$ with $\lambda_k-\lambda_{k,0} \in \ell^2(k)$ be given. Then we have{\samepage
\begin{gather}\label{lem:asymp-F2}
	|F(t,\lambda_k)-F_0(t,\lambda_k)| \in \ell^2(k)
\end{gather}
uniformly for $t\in[0,T]$.}
\end{enumerate}
\end{Lemma}

\begin{proof}For $t_0 \in [0,T]$ we consider ${q}_{t_0} \in L^2([0,T])$ defined by
\begin{gather*}
 {q}_{t_0}(t) = \begin{cases} q(t) & \text{for $0\leq t \leq t_0$}, \\ 0 & \text{for $t_0 < t \leq T$}. \end{cases}
\end{gather*}
We denote the flat $\mathfrak{sl}(2,\bbC)$-connection of equation~\eqref{eq:alpha_q} corresponding to ${q}_{t_0}$ by ${\alpha}_{t_0} = \tfrac12(\lambda\varepsilon + q_{t_0}\varepsilon_+ + \bar{q}_{t_0}\varepsilon_-)$. Then the solution of $\tfrac{{\rm d}\ }{{\rm d}t} F_{t_0} = {F}_{t_0}{\alpha}_{t_0}$ with $F_{t_0}(0,\lambda)=\one$ is given by
\begin{gather}\label{lem:asymp-F:eq:wtF}
{F}_{t_0}(t,\lambda) = \begin{cases} F(t,\lambda) & \text{for $0\leq t \leq t_0$}, \\ F(t_0,\lambda) F_0(t-t_0,\lambda) & \text{for $t_0 < t \leq T$} . \end{cases}
\end{gather}

For the proof of (i), let $\varepsilon>0$ be given. Because $\{{q}_{t_0}|t_0\in[0,T]\}$ is a relatively compact subset of $L^2([0,T])$, there exists by Theorem~\ref{thm:asymp}(1) $R>0$ so that for any $\lambda\in\bbC$ with $|\lambda|\geq R$ and any $t_0\in [0,T]$ we have
\begin{gather*}
 |F_{t_0}(T,\lambda)-F_0(T,\lambda)| \leq \varepsilon |F_0(T,\lambda)| .
\end{gather*}
We then have by equations~\eqref{lem:asymp-F:eq:wtF} and \eqref{thm:asymp:eq:F0-homo}
\begin{align*}
|F(t_0,\lambda)-F_0(t_0,\lambda)|
& = \big|F_{t_0}(T,\lambda) F_0(T-t_0,\lambda)^{-1} - F_0(T,\lambda) F_0(T-t_0,\lambda)^{-1}\big| \\
& \leq|F_{t_0}(T,\lambda) - F_0(T,\lambda)| \big|F_0(T-t_0,\lambda)^{-1}\big| \\
& \leq \varepsilon |F_0(T,\lambda)| \big|F_0(T-t_0,\lambda)^{-1}\big| = \varepsilon |F_0(t_0,\lambda)| \leq \varepsilon |M_0(\lambda)| .
\end{align*}

For the proof of (ii), we note that by Corollary~\ref{cor:asymp} applied to the potential ${q}_{t_0}$ we have
\begin{gather*}
|{F}_{t_0}(T,\lambda_k)-F_0(T,\lambda_k)| \in \ell^2(k) ,
\end{gather*}
and this estimate again holds uniformly for $t_0\in [0,T]$. By multiplying this estimate with $F_0(T-t_0,\lambda_k)^{-1}$ and noting that $\big|F_0(T-t_0,\lambda_k)^{-1}\big|$ is bounded with respect to both $t_0$ and $k$ we obtain
\begin{gather*}
 \big|{F}_{t_0}(T,\lambda_k) F_0(T-t_0,\lambda_k)^{-1} - F_0(t_0,\lambda_k)\big| \in \ell^2(k) .
\end{gather*}
From this inequality \eqref{lem:asymp-F2} follows by means of equation~\eqref{lem:asymp-F:eq:wtF}.
\end{proof}

We are now interested in the variation of the monodromy $M$ at a point $q\in L^2([0,T])$. We fix $q\in L^2([0,T])$, and denote infinitesimal variations of $q$ (in other words, tangent vectors to $L^2([0,T])$ at~$q$) by $\delta q$. For a functional $f$ defined on $L^2([0,T])$ we denote the corresponding variation of $f$ by $\delta f = \tfrac{\mathrm{d}f}{\mathrm{d}q} \delta q$. For $q \in L^2([0,T])$
and $k\in \bbZ$ we denote by
\begin{gather*}
 \widehat{q}(k) := \int_0^T q(t) \exp\left( \tfrac{2\pi \mi}{T}kt \right) \mathrm{d}t
\end{gather*}
the $k$-th Fourier coefficient of $q$.

\begin{Proposition}\label{prop:asymp-deltaM}Let $q\in L^2([0,T])$ and a sequence $(\lambda_k)_{k\in \bbZ}$ with $\lambda_k - \lambda_{k,0} \in \ell^2(k)$ be given. Then there exists a sequence $(\rho_k)_{k\in \bbZ} \in \ell^2(k)$ so that we have for every variation $\delta q$ of~$q$
\begin{gather}\label{prop:asymp-deltaM:eq:estimate}
\left| 2(-1)^k\delta M(\lambda_k) - \left( \begin{smallmatrix} 0 & \widehat{\delta q}(k) \\ -\widehat{\delta \bar{q}}(-k) & 0 \end{smallmatrix} \right) \right| \leq \|\delta q\|_{L^2([0,L])} \rho_k .
\end{gather}
\end{Proposition}

\begin{proof}Because the left-hand side of~\eqref{prop:asymp-deltaM:eq:estimate} is homogeneous with respect to $\delta q$, it suffices to consider the case $\|\delta q\|_{L^2([0,T])} \leq 1$. Then we are to show that
\begin{gather}\label{prop:asymp-deltaM:eq:estimate2}
\left| \delta M(\lambda_k) - \tfrac{(-1)^k}{2}\left( \begin{smallmatrix} 0 & \widehat{\delta q}(k) \\ -\widehat{\delta \bar{q}}(-k) & 0 \end{smallmatrix} \right) \right| \leq \rho_k
\end{gather}
holds with some sequence $\rho_k \in \ell^2(k)$ which is independent of $\delta q$.

We have $\delta M(\lambda_k) = \delta M(\lambda_{k,0}) + \bigr( \delta M(\lambda_k)-\delta M(\lambda_{k,0}) \bigr)$, and will treat the two summands on the right-hand side of this equation separately to obtain the estimate of~\eqref{prop:asymp-deltaM:eq:estimate2}.

First we note that by taking the derivative of the differential equation $\tfrac{{\rm d} \ }{{\rm d}t} F = F\alpha$ with respect to $q$, we obtain the differential equation $\tfrac{{\rm d} \ }{{\rm d}t}(\delta F) = (\delta F) \alpha + F \delta \alpha$ for~$\delta F$. Therefore we have
\begin{gather}\label{prop:asymp-deltaM:eq:deltaF}
\delta F(t_1,\lambda)= \int_0^{t_1} F(t,\lambda)(\delta \alpha)(t,\lambda)F(t,\lambda)^{-1}\mathrm{d}t \cdot F(t_1,\lambda)
\end{gather}
and hence
\begin{align}
\delta M(\lambda_{k,0}) & = \int_0^{T} F(t,\lambda_{k,0})(\delta \alpha)(t,\lambda_{k,0})F(t,\lambda_{k,0})^{-1}\mathrm{d}t \cdot M(\lambda_{k,0}) \nonumber\\
& = \int_0^{T} F_0(t,\lambda_{k,0})(\delta \alpha)(t,\lambda_{k,0})F_0(t,\lambda_{k,0})^{-1}\mathrm{d}t \cdot M_0(\lambda_{k,0}) + r_k\label{prop:asymp-deltaM:eq:part1-pre}
\end{align}
with
\begin{gather*}
r_k := \int_0^{T} (F-F_0)(t,\lambda_{k,0})(\delta \alpha)(t,\lambda_{k,0})F_0(t,\lambda_{k,0})^{-1}\mathrm{d}t\cdot M_0(\lambda_{k,0}) \\
\hphantom{r_k :=}{} + \int_0^{T} F(t,\lambda_{k,0})(\delta \alpha)(t,\lambda_{k,0})\big(F^{-1}-F_0^{-1}\big)(t,\lambda_{k,0})\mathrm{d}t \cdot M_0(\lambda_{k,0}) \\
\hphantom{r_k :=}{} + \int_0^{T} F(t,\lambda_{k,0})(\delta \alpha)(t,\lambda_{k,0})F(t,\lambda_{k,0})^{-1}\mathrm{d}t \cdot (M-M_0)(\lambda_{k,0}) .
\end{gather*}
By equation~\eqref{eq:alpha_q} we have
\begin{gather}\label{prop:asymp-deltaM:eq:deltaalpha}
\delta\alpha = \tfrac12(\delta q \varepsilon_+ + \delta\bar{q}\varepsilon_-) ,
\end{gather}
therefore $\delta \alpha$ is independent of $\lambda$, $\delta \alpha \in L^2\big([0,T],\C^{2\times 2}\big)$ holds, and we have $\|\delta \alpha\|_{L^2([0,T],\C^{2\times 2})} = \|\delta q\|_{L^2([0,T])} \leq 1$. Because of Lemma~\ref{lem:asymp-F}(2) and the fact that $F(t,\lambda_{k,0})$ is bounded uniformly with respect to $t\in[0,T]$ and $k\in\bbZ$, it follows that $r_k \in \ell^2(\bbZ)$ holds. Moreover
there exists a sequence $\rho_k^{(1)} \in \ell^2(\bbZ)$ that is independent of $\delta q$ so that $|r_k| \leq \rho_k^{(1)}$ holds.

By using equations~\eqref{eq:M0}, \eqref{eq:lambdak0} and \eqref{prop:asymp-deltaM:eq:deltaalpha} as well as $M_0(\lambda_{k,0}) = (-1)^k \one$, we obtain from equation~\eqref{prop:asymp-deltaM:eq:part1-pre}
\begin{align}
\delta M(\lambda_{k,0}) & = \frac{(-1)^k}{2} \int_0^{T} \exp\left( \frac{\pi kt}{T} \varepsilon \right) \begin{pmatrix} 0 & \delta q(t) \\ -\overline{\delta q(t)} & 0 \end{pmatrix} \exp\left( -\frac{\pi kt}{T} \varepsilon \right)\mathrm{d}t + r_k \nonumber\\
& = \frac{(-1)^k}{2} \int_0^{T} \begin{pmatrix} 0 & \delta q(t)\exp\left( \tfrac{2\pi \mi}{T}kt \right) \\ -\overline{\delta q(t)}\exp\left( -\tfrac{2\pi \mi}{T}kt \right) & 0 \end{pmatrix}\mathrm{d}t + r_k \nonumber\\
& = \frac{(-1)^k}{2} \begin{pmatrix} 0 & \widehat{\delta q}(k) \\ -\widehat{\delta \bar{q}}(-k) & 0 \end{pmatrix} + r_k .\label{prop:asymp-deltaM:eq:part1}
\end{align}
Moreover we have
\begin{gather}\label{prop:asymp-deltaM:eq:part2-pre}
\delta M(\lambda_k)-\delta M(\lambda_{k,0}) = \int_{\lambda_{k,0}}^{\lambda_k} (\delta M)'(\lambda)\mathrm{d}\lambda ,
\end{gather}
where $'$ denotes differentiation with respect to $\lambda$, and by taking the derivative with respect to $\lambda$ of equation~\eqref{prop:asymp-deltaM:eq:deltaF} and setting $t=T$ we obtain
\begin{gather*}
(\delta M)'(\lambda) = \int_0^T F'(t,\lambda) (\delta \alpha)(t,\lambda)F(t,\lambda)^{-1}\mathrm{d}t \cdot M(\lambda) \\
\hphantom{(\delta M)'(\lambda) =}{} + \int_0^T F(t,\lambda) (\delta \alpha)(t,\lambda)(F^{-1})'(t,\lambda)\mathrm{d}t \cdot M(\lambda) \\
\hphantom{(\delta M)'(\lambda) =}{} + \int_0^T F(t,\lambda) (\delta \alpha)(t,\lambda)F(t,\lambda)^{-1}\mathrm{d}t \cdot M'(\lambda) .
\end{gather*}
It follows similarly as in the proof of Corollary~\ref{cor:asymp-M'} that $F$ and $F^{-1}$ are bounded uniformly for $\lambda \in \bigcup_{k\in \bbZ} [\lambda_{k,0},\lambda_k]$ and $t\in [0,T]$; because moreover $\delta \alpha \in L^2\big([0,T],\C^{2\times 2}\big)$ holds, where $\delta \alpha$ does not depend on $\lambda$, the preceding equation shows that $(\delta M)'(\lambda)$ is bounded on~$\lambda \in \bigcup_{k\in \bbZ} [\lambda_{k,0},\lambda_k]$. Thus it follows from equation~\eqref{prop:asymp-deltaM:eq:part2-pre} that
\begin{gather}\label{prop:asymp-deltaM:eq:part2}
|\delta M(\lambda_k)-\delta M(\lambda_{k,0})| \leq C |\lambda_k-\lambda_{k,0}| =: \rho_k^{(2)}
\end{gather}
holds with some $C>0$. Because of $\lambda_k-\lambda_{k,0} \in \ell^2(k)$, we have $\rho_k^{(2)} \in \ell^2(k)$.

It follows from \eqref{prop:asymp-deltaM:eq:part1} and~\eqref{prop:asymp-deltaM:eq:part2} that the estimate \eqref{prop:asymp-deltaM:eq:estimate2} holds with $\rho_k := \rho_k^{(1)} + \rho_k^{(2)} \in \ell^2(k)$.
\end{proof}

\section{The perturbed Fourier map}\label{Se:perturbed-fourier-map}

We now consider the map
\begin{gather*}
 \Phi\colon \ L^2([0,T]) \to \ell^2(\bbZ), q \mapsto (z_k)
\end{gather*}
that associates to each potential $q\in L^2([0,T])$ the corresponding perturbed Fourier coefficients (Definition~\ref{def:perturbedfourier}). It is clear that $\Phi$ is smooth in the ``weak'' sense that each component map $L^2([0,T])\to \bbC$, $q \mapsto z_k$ is smooth. We will see in this section that $\Phi$ is in fact smooth as a map of Banach spaces, and that by restricting it to suitable affine subspaces of finite co-dimension, we obtain a local diffeomorphism.

\begin{Lemma}\label{lem:delta-zk} Let $q\in L^2([0,T])$ be given. Then there exists $\rho_k \in \ell^2(k)$ so that for any variation $\delta q \in L^2([0,T])$ of $q$ with $\|\delta q\|_{L^2([0,T])} \leq 1$, the corresponding variation $\delta z_k$ at $q$ in the direction $\delta q$ of the perturbed Fourier coefficient satisfies
\begin{gather}\label{lem:delta-zk:eq:claim}
\big| \delta z_k - \widehat{\delta q}(k) \big| \leq \rho_k .
\end{gather}
\end{Lemma}

\begin{proof}Let $(\lambda_k)$ be as in Lemma~\ref{lem:lambdak}. Because the $\lambda_k$ are characterised by the equation $(a-d)(\lambda_k)=0$, we have
\begin{gather*}
 \delta\lambda_k = - \frac{1}{(a-d)'(\lambda_k)} {\delta(a-d)(\lambda_k)} .
\end{gather*}
We have $(a_0-d_0)'(\lambda_k) = \mi T\cosh\big(\tfrac{\mi T}{2}\lambda_k\big) = 2(-1)^k\mi T + \ell^2(k)$, and therefore Corollary~\ref{cor:asymp-M'} shows that $(a-d)'(\lambda_k)$ is bounded away from zero,
whence it follows that $\tfrac{1}{(a-d)'(\lambda_k)}$ is bounded with respect to $k\in \bbZ$. Moreover it follows from Proposition~\ref{prop:asymp-deltaM} that there exists $\rho_k^{(1)} \in \ell^2(k)$
(independent of $\delta q$) so that $|\delta(a-d)(\lambda_k)| \leq \rho_k^{(1)}$ holds. Thus we see that there exists $\rho_k^{(2)} \in \ell^2(k)$ (again independent of $\delta q$) so that{\samepage
\begin{gather}\label{lem:delta-zk:eq:lambdak}
|\delta \lambda_k| \leq \rho_k^{(2)}
\end{gather}
holds for all $k\in \bbZ$.}

By definition we have $z_k = 2(-1)^kb(\lambda_k)$, and therefore
\begin{gather}\label{lem:delta-zk:eq:zk}
\delta z_k = 2(-1)^k\bigr(\delta b(\lambda_k) + b'(\lambda_k) \delta \lambda_k\bigr) .
\end{gather}
Again by Proposition~\ref{prop:asymp-deltaM} there exists $\rho_k^{(3)} \in \ell^2(k)$ (once again independent of $\delta q$) so that
\begin{gather}\label{lem:delta-zk:eq:zk-1}
\big| 2(-1)^k\delta b(\lambda_k) - \widehat{\delta q}(k) \big| \leq \rho_k^{(3)}
\end{gather}
holds. Moreover, $b'(\lambda_k)$ is bounded because of Corollary~\ref{cor:asymp-M'}, and thus by \eqref{lem:delta-zk:eq:lambdak} we have
\begin{gather}\label{lem:delta-zk:eq:zk-2}
\big|2(-1)^kb'(\lambda_k) \delta \lambda_k\big| \leq C \rho_k^{(2)}
\end{gather}
with some $C>0$. Plugging \eqref{lem:delta-zk:eq:zk-1}, \eqref{lem:delta-zk:eq:zk-2} into~\eqref{lem:delta-zk:eq:zk} yields~\eqref{lem:delta-zk:eq:claim} with $\rho_k := \rho_k^{(3)} + C \rho_k^{(2)} \in \ell^2(k)$.
\end{proof}

\begin{Proposition}\label{prop:Phi-diffeo}Let $q\in L^2([0,T])$ be given. Then there exists $N \in \bbN$ so that the map
\begin{gather*}
 \Phi_N\colon \ L^2([0,T])_{q,N} \to \ell^2(|k|> N), \qquad q_1 \mapsto (z_k(q_1))_{|k|> N}
\end{gather*}
is a local diffeomorphism near $q$, where
\begin{gather}\label{prop:Phi-diffeo:eq:L2qN}
L^2([0,T])_{q,N} := \big\{ q_1\in L^2([0,T]) \,|\, \widehat{q_1}(k)=\widehat{q}(k) \text{ for all $|k|\leq N$} \big\} .
\end{gather}
\end{Proposition}

\begin{proof}Let $\rho_k \in \ell^2(k)$ be as in Lemma~\ref{lem:delta-zk}, and then choose $N\in \bbN$ so large that the $\ell^2$-norm of the end piece sequence $(\rho_k)_{|k|>N}$ satisfies
\begin{gather*}%\label{prop:Phi-diffeo:eq:N}
\big\| (\rho_k)_{|k|>N} \big\|_{\ell^2} < \tfrac{1}{2} .
\end{gather*}
It then follows from Lemma~\ref{lem:delta-zk} that the ``weak'' derivative
\begin{gather*}
 \Phi_N'(q)\colon \ L^2([0,T])_{0,N} \to \ell^2(|k|>N), \qquad \delta q \mapsto (\delta z_k)_{|k|>N}
\end{gather*}
satisfies $\|\Phi_N'(q)-\mathfrak{F}_N\|<\tfrac12$, where $\mathfrak{F}_N\colon L^2([0,T])_{0,N} \to \ell^2(|k|>N)$, $q_1 \mapsto (\widehat{q_1}(k))_{|k|>N}$ is the ordinary Fourier transform. Because $\mathfrak{F}_N$ is an isometry, it follows that~$\Phi_N$ is differentiable at~$q$ as a map of Banach spaces, and that the Banach space homomorphism~$\Phi'(q)$ is invertible. Hence~$\Phi_N$ is a local diffeomorphism near~$q$ by the inverse function theorem.
\end{proof}

In the following corollary, we apply Proposition~\ref{prop:Phi-diffeo} to actually construct finite gap potentials that are $L^2$-close to a given potential $q\in L^2([0,T])$.

\begin{Corollary}\label{cor:dense}The set of finite gap potentials is dense in $L^2([0,T])$.
\end{Corollary}

\begin{proof}Let $q\in L^2([0,T])$ be given. By Proposition~\ref{prop:Phi-diffeo} there exist $N\in\bbN$, neighborhoods $V$ of $q$ in $L^2([0,T])_{q,N}$ and $W$ of $(z_k)_{|k|>N} := \Phi_N(q)$ in $\ell^2(|k|>N)$,
so that $\Phi_N|V\colon V \to W$ is a~diffeomorphism. For each $n\in \bbN$ with $n\geq N$, we define a sequence $\big(z^{(n)}_k\big)_{|k|>N}$ in $\ell^2(|k|>N)$ by
\begin{gather*}
 z^{(n)}_k := \begin{cases} z_k & \text{if $|k|\leq n$}, \\ 0 & \text{if $|k|>n$} . \end{cases}
\end{gather*}
In $\ell^2(|k|>N)$, the sequence $\big(z^{(n)}_k\big)_{|k|>N}$ then converges for $n\to\infty$ to $(z_k)_{|k|>N}$, therefore there exists $N_1\geq N$ so that we have $\big(z^{(n)}_k\big)_{|k|>N} \in W$ for all $n>N_1$. For such $n$ we put $q_n := (\Phi_N|V)^{-1}\big((z^{(n)}_k)_{|k|>N}\big)$. Because only finitely many of the perturbed Fourier coefficients of~$q_n$ are non-zero, $q_n$ is a finite gap
potential, and because $\Phi_N|V\colon V \to W$ is a diffeomorphism, $(q_n)_{n>N_1}$ converges to $q$ in $L^2([0,T])$.
\end{proof}

It is well-known that finite gap potentials always extend to $\bbR$ and are smooth. It is a~remarkable fact that the finite gap potentials constructed in the proof of the preceding corollary (seen as smooth functions on~$\bbR$) are always periodic with period~$T$. This is shown in the following proposition:

\begin{Proposition}\label{prop:finitegap-periodic} Let $q\in L^2([0,T])$ be a finite gap potential. Then $q$ extends to a finite gap potential on $\bbR$ which is smooth and is periodic with period $T$.
\end{Proposition}

\begin{proof}It is well-known that $q$ extends to a smooth, finite gap potential on~$\R$, which we will also denote by~$q$. This follows from the explicit description of finite gap potentials $q$ in terms of the Riemann theta function, see \cite[equation~(36)]{calini-ivey2005}, and see also~\cite{algebro}. We need to show that $q$ is periodic with period $T$.

We will consider a normalized eigenfunction of the monodromy, that is a meromorphic section of the eigenline bundle $\Lambda$ on $\widehat{\Sigma}$ (see the end of Section~\ref{Se:perturbed-fourier} for the definition of $\widehat{\Sigma}$). For the monodromy $M(\lambda)$ of~$q$, it would be possible that the section has a pole at~$\infty_\pm$, which would prevent the following approach from working. To overcome this complication, we regauge $\alpha$ from the right with the constant matrix $g:= \tfrac{1}{\sqrt{2}} \left( \begin{smallmatrix} 1 & -\mi \\ 1 & \mi \end{smallmatrix} \right)$, giving $\widetilde{\alpha} := \alpha.g = g^{-1}\alpha g$. Then the solution of $\tfrac{{\rm d} \ }{{\rm d}t} \widetilde{F}=\widetilde{F}\widetilde{\alpha}$ with $\widetilde{F}(0,\lambda)=\one$ is given by $\widetilde{F} = g^{-1}F g$, and the corresponding monodromy is $\widetilde{M}=g^{-1}M g$. The eigenvalues of~$M$ and of~$\widetilde{M}$ are the same, therefore the eigenvectors of~$\widetilde{M}$ define a holomorphic line bundle $\widetilde{\Lambda}$ on the same complex curve $\widehat{\Sigma} \setminus \{\infty_\pm\}$ as before.

Writing $\widetilde{M}(\lambda) = \left( \begin{smallmatrix} \widetilde{a}(\lambda) & \widetilde{b}(\lambda) \\ \widetilde{c}(\lambda) & \widetilde{d}(\lambda) \end{smallmatrix} \right)$, the meromorphic function $\widetilde{v} = \big(1,\tfrac{\mu-\widetilde{a}}{\widetilde{b}}\big)^t$ on $\widehat{\Sigma}$ is a meromorphic section in the eigenline bundle $\widetilde{\Lambda}$. We consider the meromorphic function $\Psi := \widetilde{F} \widetilde{v} \colon \bbR \times \big(\widehat{\Sigma} \setminus \{\infty_\pm\}\big) \to \bbC^2$.

By Lemma~\ref{lem:asymp-F}(1), $F$ is asymptotically close to the extended frame $F_0$ of the vacuum. Therefore $\widetilde{F}$, $\widetilde{v}$ and $\Psi$ are also asymptotically close to the corresponding quantities $\widetilde{F}_0$, $\widetilde{v}_0$ and $\Psi_0$ for the vacuum. We have
\begin{gather*}
 \widetilde{F}_0(t,\lambda) = g^{-1}F_0(t,\lambda)g = \begin{pmatrix} \cos\left( \tfrac{1}{2}\lambda t \right) & \sin\left( \tfrac{1}{2}\lambda t \right) \\ -\sin\left( \tfrac{1}{2}\lambda t \right) & \cos\left( \tfrac{1}{2}\lambda t \right) \end{pmatrix},
 \qquad
 \widetilde{v}_0(\lambda) = \left( \begin{smallmatrix} 1 \\ (\mu_0-\widetilde{a}_0)/\widetilde{b}_0 \end{smallmatrix} \right) = \left( \begin{smallmatrix} 1 \\ \pm \mi \end{smallmatrix} \right) ,
\end{gather*}
and therefore
\begin{gather*}
 \Psi_0 = \widetilde{F}_0 \widetilde{v}_0 = \exp\left( \pm \tfrac{\mi}{2}\lambda t \right) \left( \begin{smallmatrix} 1 \\ \pm \mi \end{smallmatrix} \right) .
\end{gather*}
This shows that the function $\Psi$ is uniquely characterised as being the Baker--Akhiezer function (see \cite[Definition~8.7]{klss2016}) on $\widehat{\Sigma}$ with the following data: The ``marked points'' on $\widehat{\Sigma}$ are the two points $q_1 = \infty_+$ and $q_2 = \infty_-$ above $\infty \in \CPone$
(where we suppose that $\infty_\pm$ is labelled such that the sign corresponds to the sign in $\widetilde{v}_0 = (1,\pm \mi)^t$) with the Mittag--Leffler distribution at $(q_1,q_2)$ given by $h = \big(\tfrac{1}{4\pi}\lambda,-\tfrac{1}{4\pi}\lambda\big)$; away from the marked points, $\Psi$ is a holomorphic section of the locally free generalised divisor $\mathcal{S}$ corresponding to the holomorphic line bundle $\widetilde{\Lambda}$, which is the polar divisor of $\widetilde{v}$.

It follows that $\Psi$ is periodic with period~$T$ if and only if the family of line bundles $\mathcal{L}_h(t)$ induced by $h$ via the Krichever construction, see \cite[Section~7]{klss2016},
is periodic with this period. To show that the latter statement holds true, we will use \cite[Lemma~7.3(ii)]{klss2016}. In fact, $k := \tfrac{1}{2\pi\mi T}\ln(\mu)$ is a multi-valued function
on $\widehat{\Sigma}\setminus \{\infty_\pm\}$, whose values over a point differ by an element of $\tfrac{1}{T}\cdot \bbZ$. $k$~is again asymptotically close to the corresponding function for the vacuum, which is $k_0 := \pm\tfrac{1}{4\pi}\lambda$. This shows that $k$ is meromorphic at $\infty_\pm$, and is a solution of the Mittag--Leffler distribution given by $h$. It follows by \cite[Lemma~7.3(ii)]{klss2016} that $\mathcal{L}_h(t)$ is periodic with period $T$, and therefore, $\Psi$ is periodic with this period.

Because $\widetilde{v}$ is a section of the eigenline bundle of $\widetilde{F}$, it follows from the $T$-periodicity of $\Psi = \widetilde{F} \widetilde{v}$ that $\widetilde{F}$ is also $T$-periodic. Therefore $F = g\widetilde{F}g^{-1}$ is $T$-periodic, hence $\alpha = \big(\tfrac{{\rm d} \ }{{\rm d}t}F\big) F^{-1}$ is $T$-periodic, and thus $q$ is $T$-periodic.
\end{proof}

The preceding results do not take the closing conditions for the curve corresponding to the potential~$q$ into account. To obtain an analogous result for potentials which satisfy a closing condition as in Remark~\ref{th:monodromy-remarks}(vi), we study the variation of the closing condition in the following section.

\section{Variations of the closing conditions}\label{Se:closing}

Let $q\in L^2([0,T])$ be a periodic potential with extended frame $F$ and monodromy~$M$. The corresponding spectral curve $\Sigma$ is given by equation~\eqref{prop:finitegap-periodic:spectralcurve}. It is hyperelliptic, with the hyperelliptic involution being given by $\sigma\colon (\lambda,\mu) \mapsto \big(\lambda,\mu^{-1}\big)$. The eigenbundle~$\Lambda$ of~$M$ is a holomorphic line bundle on a certain partial desingularisation $\widehat{\Sigma}$ of~$\Sigma$, as was described in the proof of Proposition~\ref{prop:finitegap-periodic}. The eigenline bundle $\Lambda^t$ of the transpose $M^t$ is also a holomorphic line bundle on~$\widehat{\Sigma}$. We denote non-trivial holomorphic sections of $\Lambda$ and $\Lambda^t$ by $v$ and~$w$, respectively; they are maps $v,w\colon \widehat{\Sigma}\to \bbC^2$. In terms of these sections, the projection operator onto the eigenline bundle of~$M(\lambda)$ is given by
\begin{gather*}
 P := \frac{v w^t}{w^t v} ;
\end{gather*}
note that $P$ exists as a holomorphic operator even at those points of $\widehat{\Sigma}$ where $w^t v=0$ holds (this is the case whenever $\mu=\pm 1$). The following lemma shows that the eigenvalue~$\mu$ as a~holomorphic function on $\widehat{\Sigma}$ can be recovered from $v$ and $w$ via the projector $P$.

\begin{Lemma}\label{lem:A-mu}\quad
\begin{enumerate}\itemsep=0pt
\item[{\rm{(i)}}] For any linear operator $A$ on $\bbC^2$ we have $\operatorname{tr}(P A) = \tfrac{w^t Av}{w^t v}$.
\item[{\rm{(ii)}}] $\mu = \operatorname{tr}(P M)$.
\end{enumerate}
\end{Lemma}

\begin{proof}Because both sides of both equations in the lemma are holomorphic, it suffices to show the equations for those values of $\lambda$ where $w^t v \neq 0$ holds. The kernel of $P$ is spanned by $u := v \circ \sigma^{-1}$. In fact, writing $M=\left( \begin{smallmatrix} a & b \\ c & d \end{smallmatrix} \right)$, $v=(b,\mu-a)^t$ and $w=(c,\mu-a)^t$ are holomorphic sections of the eigenline bundle of $M(\lambda)$ and $M(\lambda)^t$, respectively. Moreover, we have the equation
\begin{gather*}
 (\mu-a)(\mu-d) = bc
\end{gather*}
and therefore
\begin{align*}
 w^t u & = (c,\mu-a)\cdot \big(b,\mu^{-1}-a\big)^t = cb+(\mu-a)\big(\mu^{-1}-a\big) \\
 & = cb+(\mu-a)((\Delta-\mu)-(\Delta-d))= cb-(\mu-a)(\mu-d)=0 ,
\end{align*}
whence $Pu = 0$ follows. Because of $w^t v\neq 0$, $(v,u)$ is linear independent, and hence a basis of~$\C^2$.

Now suppose that $A$ is a linear operator on $\bbC^2$. We then have
\begin{gather*}
 PAv = \frac{v w^t Av}{w^t v} = \left(\frac{w^t Av}{w^t v}\right) v \qquad\text{and}\qquad PAu \in \bbC v ,
\end{gather*}
which implies (i). By applying (i) with $A=M$ we obtain (ii):
\begin{gather*}
 \operatorname{tr}(P M) = \frac{w^t Mv}{w^t v} = \frac{w^t \mu v}{w^t v} = \mu .\tag*{\qed}
\end{gather*}\renewcommand{\qed}{}
\end{proof}

\begin{Lemma}\label{lem:closing} Let $q\in L^2([0,T])$ with extended frame $F$ and monodromy $M$. Let $\mu$ be the eigenvalue function. For $\theta\in\R$, the closing conditions of Remark~{\rm \ref{th:monodromy-remarks}(vi)} are then equivalent to the following conditions: $(i)$ for $\bbS^3$: $\mu(1+\theta)=\mu(-1+\theta)=\pm 1$ and $(ii)$ for $\bbR^3$: $\mu(\theta)=\pm 1$, $\mu'(\theta)=0$.
\end{Lemma}
\begin{proof}
For $\lambda\in \R$ we have $M(\lambda) \in \mathrm{SU}_2$ by Remark~\ref{th:monodromy-remarks}(ii), and therefore $M(\lambda )$ is semisimple for $\lambda \in \R$.
This means that for $\lambda \in \R$, $\mu(\lambda)=\pm 1$ is equivalent to $M(\lambda)=\pm \one$. This proves (i). Diagonalizing $M = T \left( \begin{smallmatrix} \mu & 0 \\ 0 & 1/\mu \end{smallmatrix} \right)T^{-1}$ and differentiating with respect to $\lambda$ gives
$M'= \mu' T \bigl( \begin{smallmatrix} 1 & 0 \\ 0 & -1/\mu^2 \end{smallmatrix} \bigr) T^{-1} + [T'T^{-1},M] $. This shows that if $M(\theta)=\pm \one$ holds, then $M'(\theta)=0$ is equivalent
to $\mu'(0)=0$, proving (ii).
\end{proof}

\begin{Lemma} \label{lem:deltamu}Let $q\in L^2([0,T])$. We use the notations of the preceding lemma, fix a~variation $\delta q \in L^2([0,T])$ of $q$, and denote the associated variations of quantities~$f$ depending on $q$ by $\delta f$ as in Section~{\rm \ref{Se:var-monodromy}}. For each given $\lambda\in\bbC$ and eigenvalue $\mu \in \bbC$ of $M(\lambda)$ there exist non-zero eigenvectors $v,w\in\bbC^2$ so that $M(\lambda)v=\mu v$, $M(\lambda)^t w=\mu w$, and $w^t v \neq 0$. For $x\in [0,T]$ we define
\begin{gather*}
 v(x) := F(x,\lambda)^{-1} v \qquad\text{and}\qquad w(x) := F(x,\lambda)^t w .
\end{gather*}
Then we have
\begin{gather*}
 \delta \mu = \frac{\mu}{w^t v} \int_0^T w(t)^t\delta\alpha(t)v(t)\mathrm{d}t \qquad\text{where}\qquad \delta\alpha = \tfrac12 (\delta q \varepsilon_+ + \delta\bar q\varepsilon_-) .
\end{gather*}
\end{Lemma}

\begin{proof}By variation of the equation of Lemma~\ref{lem:A-mu}(ii) we obtain
\begin{gather}\label{lem:deltamu:eq:deltamu-1}
\delta \mu = \delta\bigr( \operatorname{tr}(P M(\lambda))\bigr) = \operatorname{tr}(P \delta M(\lambda)) + \operatorname{tr}(\delta P M(\lambda)) .
\end{gather}
$P$ being a projection operator, we have $P = P^2$ and therefore $\delta P = P\cdot \delta P + \delta P \cdot P$. The latter equation implies $P\cdot (\delta P)v = 0$ hence $(\delta P) v \in \C u$ (where we again put $u := v \circ \sigma$ as in the proof of Lemma~\ref{lem:A-mu}), and also $(\delta P) u = P (\delta P)u$ hence $(\delta P) u \in \C v$. Because both $v$ and $u$ are eigenvalues of $M(\lambda)$, it follows that $\delta P \cdot M(\lambda)v \in \C u$ and $\delta P \cdot M(\lambda)u \in \C v$ holds, whence
$\operatorname{tr}(\delta P \cdot M(\lambda))=0$ follows. Thus it follows from \eqref{lem:deltamu:eq:deltamu-1} that
\begin{gather} \label{lem:deltamu:eq:deltamu-2}
\delta \mu = \operatorname{tr}(P\cdot \delta M(\lambda)) .
\end{gather}

We now calculate $\delta M(\lambda)$. For this purpose we note that as consequence of the initial value problem for $F$: $\tfrac{{\rm d} \ }{{\rm d}t}F = F \alpha$, $F(0,\lambda)=\one$, $\delta F$ is characterised by the initial value problem
\begin{gather*}
 \tfrac{{\rm d} \ }{{\rm d}t} \delta F = \delta F \alpha + F \delta \alpha \qquad\text{with}\quad \delta F(0,\lambda)=0 ,
\end{gather*}
which has the unique solution
\begin{gather*}
 \delta F(x,\lambda) = \int_0^x F(t,\lambda)\delta\alpha(t)F(t,\lambda)^{-1}\mathrm{d}t \cdot F(x,\lambda) .
\end{gather*}
Therefore
\begin{gather*}
 \delta M(\lambda) = \int_0^T F(t,\lambda)\delta\alpha(t)F(t,\lambda)^{-1}\mathrm{d}t \cdot M(\lambda)
\end{gather*}
holds. We therefore obtain from equation~\eqref{lem:deltamu:eq:deltamu-2} and Lemma~\ref{lem:A-mu}(i) (applied with $A=\delta M(\lambda)$)
\begin{align*}
\delta \mu & = \operatorname{tr}(P \delta M(\lambda)) = \frac{w^t \delta M(\lambda)v}{w^t v}
= \frac{1}{w^t v} \left( w^t \int_0^T F(t,\lambda)\delta\alpha(t)F(t,\lambda)^{-1}\mathrm{d}t \cdot M(\lambda)v \right) \\
& = \frac{\mu}{w^t v} \int_0^T w^t F(t,\lambda)\delta\alpha(t)F(t,\lambda)^{-1}v \mathrm{d}t
= \frac{\mu}{w^t v} \int_0^T w(t)^t \delta\alpha(t)v(t) \mathrm{d}t .\tag*{\qed}
\end{align*}\renewcommand{\qed}{}
\end{proof}

\begin{Lemma}\label{lem:closing-linindep}Suppose that $q\in L^2([0,T])$ is not of the form $q(x)=a{\rm e}^{cx}$ with constants $a,c\in\bbC$, and that it satisfies either of the two closing conditions in Lemma~{\rm \ref{lem:closing}} for some $\theta\in\R$.
\begin{enumerate}\itemsep=0pt
\item[{\rm{(i)}}] If $q$ satisfies the closing condition for $\bbS^3$ $($Lemma~{\rm \ref{lem:closing}(i))}, there exist two variations $\delta_1 q,\delta_2 q$ $\in L^2([0,T])$ of $q$ so that the matrix
\begin{gather}\label{lem:closing-linindep:S3:matrix}
\begin{pmatrix} \delta_1 \mu(1+\theta) & \delta_2 \mu(1+\theta) \\ \delta_1 \mu(-1+\theta) & \delta_2 \mu(-1+\theta) \end{pmatrix}
\end{gather}
has maximal rank.
\item[{\rm{(ii)}}] If $q$ satisfies the closing condition for $\bbR^3$ $($Lemma~{\rm \ref{lem:closing}(ii))} there exist two variations $\delta_1 q,\delta_2 q$ $\in L^2([0,T])$ of $q$ so that the matrix
\begin{gather*}%\label{lem:closing-linindep:R3:matrix}
\begin{pmatrix} \delta_1 \mu(\theta) & \delta_2 \mu(\theta) \\ \delta_1 \mu'(\theta) & \delta_2 \mu'(\theta) \end{pmatrix}
\end{gather*}
has maximal rank.
\end{enumerate}

In both cases the variations $\delta_1 q$ and $\delta_2 q$ can be chosen so that only finitely many of their Fourier coefficients are non-zero.
\end{Lemma}

\begin{proof}The proofs for the two cases follow the same general pathway, but differ in their details due to the different closing conditions in each case. In each case, if the required variations $\delta_1 q$ and $\delta_2 q$ exist at all, then they can be chosen with only finitely many non-zero Fourier coefficients. This is true because the set of $L^2$-functions with finitely many non-zero Fourier coefficients is dense in $L^2([0,T])$ and the condition required of the $\delta_k q$ is open.

(i) Let us assume to the contrary that no such variations $\delta_1 q$ and $\delta_2 q$ exist. This means that the two linear forms
\begin{gather*}
 L^2([0,T]) \to \C, \qquad \delta q \mapsto \delta \mu(1+\theta) \qquad\text{and}\qquad L^2([0,T]) \to \C,\qquad \delta q \mapsto \delta \mu(-1+\theta)
\end{gather*}
are linear dependent over $\C$, which because of Lemma~\ref{lem:deltamu} means that the two linear forms on~$L^2([0,T])$
\begin{gather*}
 \delta q \mapsto 2 \int_0^T w(t)^t \delta\alpha(t)v(t)\mathrm{d}t = \int_0^T \bigr( v_1(t)w_2(t)\delta q(t)-v_2(t)w_1(t)\delta \overline{q}(t)\bigr)\mathrm{d}t
\end{gather*}
for $\lambda=1+\theta$ and $\lambda=-1+\theta$ are linear dependent (where $v(t)$ and $w(t)$ are defined as in Lemma~\ref{lem:deltamu}). Hence there exist constants $s_1,s_{-1}\in \C$, not both zero, so that
\begin{gather}
s_1(v_1w_2)\bigr|_{\lambda=1+\theta} + s_{-1}(v_1w_2)\bigr|_{\lambda=-1+\theta} = 0,\nonumber\\ s_1(v_2w_1)\bigr|_{\lambda=1+\theta} + s_{-1}(v_2w_1)\bigr|_{\lambda=-1+\theta} = 0\label{lem:closing-linindep:eq:S3:s+-1a}
\end{gather}
holds. We may suppose without loss of generality that $s_1\neq 0$ holds.

We define the $(2\times 2)$-matrix-valued function $B(x) := v(x)\cdot w(x)^t$ for given $\lambda\in \C$. Then the differential equation $\tfrac{{\rm d} \ }{{\rm d}t}F = F \alpha$ implies $\tfrac{{\rm d} \ }{{\rm d}t}B = [B,\alpha]$. By explicitly expressing the latter differential equation in matrix components, we obtain for the three functions
\begin{gather}\label{lem:closing-linindep:eq:S3:varphi-def}
\varphi_1 := v_1w_2 ,\qquad \varphi_2 := v_2w_1 \qquad\text{and}\qquad \varphi_3 := \tfrac12(v_1w_1-v_2w_2)
\end{gather}
the following system of differential equations (where the functions $\varphi_k$ are at least once differentiable in the Sobolev sense):
\begin{gather}\label{lem:closing-linindep:dgl1}
\dot{\varphi}_1 = -\mi \lambda\varphi_1 + q\varphi_3, \\
\label{lem:closing-linindep:dgl2}
\dot{\varphi}_2 = \mi \lambda\varphi_2 + \overline{q}\varphi_3, \\
\label{lem:closing-linindep:dgl3}
\dot{\varphi}_3 = -\tfrac12\overline{q}\varphi_1-\tfrac12q\varphi_2 ,
\end{gather}
where we abbreviate $\dot{\varphi}_k := \tfrac{{\rm d} \ }{{\rm d}t} \varphi_k$.

We will use this system of differential equations to obtain a contradiction from~\eqref{lem:closing-linindep:eq:S3:s+-1a}, which in the present notation reads
\begin{gather}\label{lem:closing-linindep:eq:S3:s+-1b}
s_1\varphi_k \bigr|_{\lambda=1+\theta} + s_{-1}\varphi_k \bigr|_{\lambda=-1+\theta} = 0 \qquad\text{for}\quad k\in \{1,2\} .
\end{gather}
From equations~\eqref{lem:closing-linindep:eq:S3:s+-1b} and \eqref{lem:closing-linindep:dgl3} we obtain
\begin{gather*}
 s_1\dot{\varphi}_3 \bigr|_{\lambda=1+\theta} + s_{-1}\dot{\varphi}_3 \bigr|_{\lambda=-1+\theta} = 0
\end{gather*}
and therefore there exists a constant $r \in \bbC$ with
\begin{gather}
\label{lem:closing-linindep:eq:S3:eq5}
s_1{\varphi}_3 \bigr|_{\lambda=1+\theta} + s_{-1}{\varphi}_3 \bigr|_{\lambda=-1+\theta} = r .
\end{gather}
On the other hand, by differentiating \eqref{lem:closing-linindep:eq:S3:s+-1b} we obtain
\begin{gather}
\label{lem:closing-linindep:eq:S3:eq6}
s_1\dot{\varphi}_k \bigr|_{\lambda=1+\theta} + s_{-1}\dot{\varphi}_k \bigr|_{\lambda=-1+\theta} = 0 \qquad\text{for}\quad k\in \{1,2\} .
\end{gather}
By plugging equation~\eqref{lem:closing-linindep:dgl1} (for both $\lambda=1+\theta$ and $\lambda=-1+\theta$) into equation~\eqref{lem:closing-linindep:eq:S3:eq6} (for $k=1$), one obtains
\begin{gather*}
 0 = -\mi\bigr( s_1(1+\theta)\varphi_1 \bigr|_{\lambda=1+\theta} + s_{-1}(-1+\theta)\varphi_1 \bigr|_{\lambda=-1+\theta} \bigr) + q \bigr( s_1{\varphi}_3 \bigr|_{\lambda=1+\theta} + s_{-1}{\varphi}_3 \bigr|_{\lambda=-1+\theta} \bigr)
\end{gather*}
and therefore by equations~\eqref{lem:closing-linindep:eq:S3:eq5} and~\eqref{lem:closing-linindep:eq:S3:s+-1b}
\begin{gather*}
 0 = -\mi\bigr( s_1(1+\theta)\varphi_1 \bigr|_{\lambda=1+\theta} - s_{1}(-1+\theta)\varphi_1 \bigr|_{\lambda=1+\theta} \bigr) + qr = -2\mi s_1 \varphi_1 \bigr|_{\lambda=1+\theta} + qr .
\end{gather*}
Because of $s_1 \neq 0$ this implies
\begin{gather}\label{lem:closing-linindep:eq:S3:eq7}
\varphi_1 \bigr|_{\lambda=1+\theta} = \frac{r}{2\mi s_1} q .
\end{gather}
Similarly, by plugging equation~\eqref{lem:closing-linindep:dgl2} into equation~\eqref{lem:closing-linindep:eq:S3:eq6} (for $k=2$), and then applying equations~\eqref{lem:closing-linindep:eq:S3:eq5} and \eqref{lem:closing-linindep:eq:S3:s+-1b}, one also obtains
\begin{gather}\label{lem:closing-linindep:eq:S3:eq8}
\varphi_2 \bigr|_{\lambda=1+\theta} = -\frac{r}{2\mi s_1} \overline{q} .
\end{gather}
By plugging equations~\eqref{lem:closing-linindep:eq:S3:eq7} and \eqref{lem:closing-linindep:eq:S3:eq8} into equation~\eqref{lem:closing-linindep:dgl3} we now obtain
\begin{gather*}
 \dot{\varphi}_3 \bigr|_{\lambda=1+\theta} = -\frac12\overline{q}\left( \frac{r}{2\mi s_1} q \right) - \frac12q\left( -\frac{r}{2\mi s_1} \overline{q} \right) = 0 .
\end{gather*}
Therefore there exists a constant $p \in \C$ so that
\begin{gather}\label{lem:closing-linindep:eq:S3:vi3-const}
\varphi_3 \bigr|_{\lambda=1+\theta} = p
\end{gather}
holds. We now note that on one hand, we obtain by plugging equations~\eqref{lem:closing-linindep:eq:S3:eq7} and \eqref{lem:closing-linindep:eq:S3:vi3-const} into equation~\eqref{lem:closing-linindep:dgl1}
\begin{gather*}
 \dot{\varphi}_1 \bigr|_{\lambda=1+\theta} = -\mi\left( \frac{r}{2\mi s_1}q \right) + qp = \left( p - \frac{r}{2s_1}\right)q ,
\end{gather*}
on the other hand, equation~\eqref{lem:closing-linindep:eq:S3:eq7} shows that~$q$ is differentiable, and by
differentiating this equation we obtain
\begin{gather*}
 \dot{\varphi}_1 \bigr|_{\lambda=1+\theta}= \frac{r}{2\mi s_1}\tfrac{{\rm d} \ }{{\rm d}t}q .
\end{gather*}
Therefore we have
\begin{gather}\label{lem:closing-linindep:eq:S3:q-dgl}
\frac{r}{2\mi s_1}\tfrac{{\rm d} \ }{{\rm d}t}q = \left( p - \frac{r}{2s_1}\right)q .
\end{gather}
If $r\neq 0$ holds, then this differential equation implies $\tfrac{{\rm d} \ }{{\rm d}t}q=cq$ with $c := \mi \big( \tfrac{2s_1p}{r}-1 \big)$, and therefore $q(x) = a{\rm e}^{cx}$ with some $a\in \C$, in contradiction to the hypothesis of the lemma that~$q$ does not have this form. If $r=0$, but $p\neq 0$ holds, then $q=0$ follows from equation~\eqref{lem:closing-linindep:eq:S3:q-dgl}, again a contradiction. Finally, if $r=p=0$ holds, then we have $\varphi_k|_{\lambda=1+\theta}=0$ for $k\in \{1,2,3\}$ by
equations~\eqref{lem:closing-linindep:eq:S3:eq7}, \eqref{lem:closing-linindep:eq:S3:eq8} and \eqref{lem:closing-linindep:eq:S3:vi3-const}. By \eqref{lem:closing-linindep:eq:S3:varphi-def}
this implies $v(x)|_{\lambda=1+\theta}=0$ or $w(x)|_{\lambda=1+\theta}=0$ for almost every $x\in [0,T]$, which in turn implies that either of the eigenvectors $v|_{\lambda=1+\theta}$ or
$w|_{\lambda=1+\theta}$ (of~$M(1)$ or~$M(1)^t$, respectively) is zero. This is also a contradiction.

(ii) Let us again assume to the contrary that no such variations $\delta_1 q$ and $\delta_2 q$ exist. This means that the two linear forms
\begin{gather*}
 L^2([0,T]) \to \C, \qquad \delta q \mapsto \delta \mu(\theta) \qquad\text{and}\qquad L^2([0,T]) \to \C, \qquad \delta q \mapsto \overline{\delta \mu'(\theta)}
\end{gather*}
are linear dependent over $\C$. By Lemma~\ref{lem:deltamu} we have
\begin{gather*}
 \delta \mu = \frac{\mu}{2w^t v} \int_0^T \bigr( v_1(t)w_2(t)\delta q(t)-v_2(t)w_1(t)\delta \overline{q}(t)\bigr)\mathrm{d}t
\end{gather*}
and therefore
\begin{gather*}
 \delta \mu' = \left( \ln\left( \frac{\mu}{2w^t v} \right)\right)' \delta \mu + \frac{\mu}{2w^t v} \int_0^T \bigr( (v_1w_2)'(t)\delta q(t)-(v_2w_1)'(t)\delta \overline{q}(t)\bigr)\mathrm{d}t .
\end{gather*}
If we again denote functions $\varphi_1$, $\varphi_2$, $\varphi_3$ as in equations~\eqref{lem:closing-linindep:eq:S3:varphi-def} (but now at $\lambda=\theta$),
our assumption therefore implies that there exist constants $s,\widetilde{s}\in \bbC$, which are not both zero, so that
\begin{gather}\label{lem:closing-linindep:eq:R3:eq7}
s \varphi_1 + \widetilde{s} \varphi_1' = s \varphi_2 + \widetilde{s} \varphi_2' = 0
\end{gather}
holds. The functions $\varphi_k$ again satisfy the system of differential equations~\eqref{lem:closing-linindep:dgl1}--\eqref{lem:closing-linindep:dgl3}. By dif\-fe\-rentiating
these equations with respect to $\lambda$, and then taking $\lambda=\theta$, we find that the six functions $\varphi_1$, $\varphi_2$, $\varphi_3$, $\varphi_1'$, $\varphi_2'$, $\varphi_3'$
are governed by the system of differential equations
\begin{gather}
\label{lem:closing-linindep:eq:R3:eq1}
\dot{\varphi}_1 = -\mi\theta\varphi_1 + q\varphi_3, \\
\label{lem:closing-linindep:eq:R3:eq2}
\dot{\varphi}_2 = \mi\theta\varphi_2 + \overline{q}\varphi_3, \\
\label{lem:closing-linindep:eq:R3:eq3}
\dot{\varphi}_3 = -\tfrac12\overline{q}\varphi_1 - \tfrac12q\varphi_2, \\
\label{lem:closing-linindep:eq:R3:eq4}
\dot{\varphi}_1' = -\mi\varphi_1 -\mi\theta\varphi_1' + q \varphi_3', \\
\label{lem:closing-linindep:eq:R3:eq5}
\dot{\varphi}_2' = \mi\varphi_2 + \mi\theta\varphi_2' + \overline{q}\varphi_3', \\
\label{lem:closing-linindep:eq:R3:eq6}
\dot{\varphi}_3' = -\tfrac12\overline{q}\varphi_1' - \tfrac12q\varphi_2'.
\end{gather}
Here we again abbreviate $\dot{\varphi}_k = \tfrac{{\rm d} \ }{{\rm d}t} \varphi_k$. We should first note that it is not possible that $\varphi_1=\varphi_2=0$ holds, because this would imply also $\varphi_3=0$ by equations~\eqref{lem:closing-linindep:eq:R3:eq3} and~\eqref{lem:closing-linindep:eq:R3:eq1}, which would be a contradiction as in~(i). This observation implies in particular that $\widetilde{s}\neq 0$ holds
because of equation~\eqref{lem:closing-linindep:eq:R3:eq7}.

From equations~\eqref{lem:closing-linindep:eq:R3:eq3}, \eqref{lem:closing-linindep:eq:R3:eq6} and \eqref{lem:closing-linindep:eq:R3:eq7} we obtain $s \dot{\varphi}_3 + \widetilde{s} \dot{\varphi}_3'=0$, and therefore there exists a constant $r\in \bbC$ so that
\begin{gather}\label{lem:closing-linindep:eq:R3:eq8}
s \varphi_3 + \widetilde{s} \varphi_3'=r
\end{gather}
holds. By differentiating equation~\eqref{lem:closing-linindep:eq:R3:eq7} with respect to~$x$ we obtain $s \dot{\varphi}_1 + \widetilde{s} \dot{\varphi}_1'=0$ and therefore by equations~\eqref{lem:closing-linindep:eq:R3:eq1} and \eqref{lem:closing-linindep:eq:R3:eq4}
\begin{align*}
 0 & = s \dot{\varphi}_1 + \widetilde{s} \dot{\varphi}_1' = s(-\mi\theta\varphi_1 + q\varphi_3) + \widetilde{s}(-\mi\varphi_1 - \mi \theta\varphi_1' + q\varphi_3') \\
 & = (-\mi\theta s - \mi\widetilde{s})\varphi_1 - \mi\theta\widetilde{s}\varphi_1' + q(s\varphi_3 + \widetilde{s}\varphi_3') ,
\end{align*}
whence it follows by equations~\eqref{lem:closing-linindep:eq:R3:eq7} and \eqref{lem:closing-linindep:eq:R3:eq8} that
\begin{gather*}
 0 = -\mi\widetilde{s}\varphi_1 + qr
\end{gather*}
and hence
\begin{gather}\label{lem:closing-linindep:eq:R3:eq10}
\varphi_1 = \frac{r}{\mi\widetilde{s}}q
\end{gather}
holds. Similarly one also obtains
\begin{gather}\label{lem:closing-linindep:eq:R3:eq11}
\varphi_2 = -\frac{r}{\mi\widetilde{s}}\overline{q} .
\end{gather}
By plugging equations~\eqref{lem:closing-linindep:eq:R3:eq10} and \eqref{lem:closing-linindep:eq:R3:eq11} into equation~\eqref{lem:closing-linindep:eq:R3:eq3}, one obtains $\dot{\varphi}_3=0$, and therefore there exists a constant $p\in \bbC$ so that $\varphi_3=p$ holds. Equation~\eqref{lem:closing-linindep:eq:R3:eq10} also shows that~$q$ is differentiable and that $\dot{\varphi}_1 = \tfrac{r}{\mi \widetilde{s}}\tfrac{{\rm d} \ }{{\rm d}t}q$ holds. By plugging these results into equation~\eqref{lem:closing-linindep:eq:R3:eq1}
we obtain
\begin{gather*}
 \frac{r}{\mi\widetilde{s}}\tfrac{{\rm d} \ }{{\rm d}t}q = pq .
\end{gather*}
If either $p\neq 0$ or $r\neq 0$ holds, then this equation implies that $q(x)=a{\rm e}^{cx}$ for some constants $a,c\in\bbC$, which contradicts the hypothesis of the lemma. If $p=r=0$ holds, then we have $\varphi_1=\varphi_2=0$ by equations~\eqref{lem:closing-linindep:eq:R3:eq10} and~\eqref{lem:closing-linindep:eq:R3:eq11}, which is also a contradiction.
\end{proof}

\begin{Theorem}\label{thm:Psi-diffeo} Suppose that $q\in L^2([0,T])$ is not of the form $q(x)=a{\rm e}^{cx}$ with constants $a,c\in\bbC$, and that it satisfies one of the two closing conditions in Lemma~{\rm \ref{lem:closing}}. For $N\in \bbN$ we define $L^2([0,T])_{q,N}$ as in equation~\eqref{prop:Phi-diffeo:eq:L2qN}.
\begin{enumerate}\itemsep=0pt
\item[{\rm{(i)}}] If $q$ satisfies the closing conditions for $\bbS^3$ with some $\theta\in\R$ there exists $N \in \bbN$ and $f_1,f_2 \in L^2([0,T])$ so that the map
\begin{align*}
 \Psi_{N,f}\colon \ L^2([0,T])_{q,N} + \bbR f_1 + \bbR f_2 & \longrightarrow \ell^2(|k|> N) \times \bbS^1 \times \bbS^1, \\
 q_1 & \longmapsto \bigr((z_k(q_1))_{|k|> N},\mu(1+\theta),\mu(-1+\theta)\bigr)
\end{align*}
is a local diffeomorphism near~$q$.
\item[{\rm{(ii)}}] If $q$ satisfies the closing conditions for $\bbR^3$ with some $\theta\in\R$ there exists $N \in \bbN$ and $f_1,f_2 \in L^2([0,T])$ so that the map
\begin{align*}
 \Psi_{N,f}\colon \ L^2([0,T])_{q,N} + \bbR f_1 + \bbR f_2 & \longrightarrow \ell^2(|k|> N) \times \bbS^1 \times \mi\bbR, \\
 q_1 & \longmapsto \bigr((z_k(q_1))_{|k|> N},\mu(\theta),\mu'(\theta)\bigr)
\end{align*}
is a local diffeomorphism near $q$.
\end{enumerate}
\end{Theorem}

\begin{proof}
For (i),
we need to show that the derivative of $\Psi_{N,f}$ at $q$
\begin{align*}
 \Psi_{N,f}'(q)\colon \ L^2([0,T])_{0,N} \oplus \R f_1 \oplus \R f_2 & \longrightarrow \ell^2(|k|>N) \oplus \mi\bbR \oplus \mi\bbR, \\
 \delta q & \longmapsto \bigr( (\delta z_k)_{|k|>N},\delta\mu(1+\theta),\delta\mu(-1+\theta)\bigr)
\end{align*}
is an isomorphism of Banach spaces. By Lemma~\ref{lem:closing-linindep}(i) there exist two variations $\delta_1 q$, $\delta_2 q$ of $q$ with only finitely many non-zero Fourier coefficients, so that the matrix \eqref{lem:closing-linindep:S3:matrix} has maximal rank. Let $f_\nu := \delta_\nu q \in L^2([0,T])$ for $\nu\in \{1,2\}$. By Proposition~\ref{prop:Phi-diffeo} there exists $N \in \bbN$ so that
\begin{gather*}
 L^2([0,T])_{0,N} \to \ell^2(|k|>N), \qquad \delta q \mapsto (\delta z_k)_{|k|>N}
\end{gather*}
is an isomorphism of Banach spaces; we can choose $N$ large enough so that additionally $\widehat{f_1}(k) = \widehat{f_2}(k)=0$ for all $k\in\bbZ$ with $|k|>N$.
Then $\Psi_{N,f}'(q)$ is an isomorphism of Banach spaces. The proof of (ii) is analogous to that of (i).
\end{proof}

\begin{Corollary}\label{cor:S3-dense}The set of finite gap potentials of $T$-periodic curves in $\bbS^3$ respectively $\bbR^3$ with some total torsion $\theta T \in \R$ is $L^2$-dense in the set of all potentials of $T$-periodic curves in $\bbS^3$ respectively $\bbR^3$ with that total torsion.
\end{Corollary}

\begin{proof} We prove the corollary for $\bbS^3$; the proof for $\bbR^3$ is analogous. Let $q\in L^2([0,T])$ be the potential of a $T$-periodic curve in $\bbS^3$. If $q$ is of the form $q(x)=a{\rm e}^{cx}$ with constants $a,c\in\bbC$, then $q$ has finite gaps, so there is nothing to show. Otherwise, Theorem~\ref{thm:Psi-diffeo}(i) shows that there exist $N \in \bbN$, $f_1,f_2 \in L^2([0,T])$, neighborhoods $V$ of $q$ in $L^2([0,T])_{q,N}+ \bbR f_1 + \bbR f_2$ and $W$ of $\bigr( (z_k)_{|k|>N},\eta,\eta\bigr) := \Psi_{N,f}(q)$ in $\ell^2(|k|>N) \times \bbS^1 \times \bbS^1$ (where $\eta := \mu(1)=\mu(-1)\in\{\pm 1\}$), so that $\Psi_{N,f}|V\colon V \to W$ is a diffeomorphism. For each $n\in \bbN$ with $n\geq N$, we define a sequence $(z^{(n)}_k)_{|k|>N}$ in $\ell^2(|k|>N)$ by
\begin{gather*}
 z^{(n)}_k := \begin{cases} z_k & \text{if $|k|\leq n$}, \\ 0 & \text{if $|k|>n$} . \end{cases}
\end{gather*}
In $\ell^2(|k|>N)$, the sequence $\bigr(\big(z^{(n)}_k\big)_{|k|>N},\eta,\eta\bigr)$ then converges for $n\to\infty$ to $\bigr((z_k)_{|k|>N},\eta,\eta\bigr)=\Psi_{N,f}(q)$, therefore there exists $N_1\geq N$ so that we have $\bigr(\big(z^{(n)}_k\big)_{|k|>N},\eta,\eta\bigr) \in W$ for all \mbox{$n>N_1$}. For such $n$ we put $q_n := (\Phi_N|V)^{-1}\bigr(\big(z^{(n)}_k\big)_{|k|>N},\eta,\eta\bigr)$. Because only finitely many of the perturbed Fourier coefficients of $q_n$ are non-zero, $q_n$ is a finite gap potential, by Lemma~\ref{lem:closing}(i), $q_n$~satisfies the closing condition for $\bbS^3$, hence corresponds to a~$T$-periodic curve in $\bbS^3$, and because $\Phi_N|V\colon V \to W$ is a~diffeomorphism, $(q_n)_{n>N_1}$ converges to~$q$ in $L^2([0,T])$.
\end{proof}

\begin{Theorem}\label{T:curvedense}The set of closed finite gap curves in $\bbS^3$ respectively $\bbR^3$ with respect to the period~$T$ is $W^{2,2}$-dense in the Sobolev space of all closed $W^{2,2}$-curves of length~$T$ in~$\bbS^3$ respectively~$\bbR^3$. Moreover, near any closed curve $\gamma$ in $\bbS^3$ or $\bbR^3$ there are closed finite gap curves with the same total torsion as~$\gamma$.
\end{Theorem}

\begin{proof} Let $\bbE^3 \in \big\{\bbS^3,\bbR^3\big\}$, and let $\gamma\colon [0,T]\to \bbE^3$ be a $W^{2,2}$-closed $W^{2,2}$-curve. We let $\theta T$ be the total torsion of $\gamma$, and consider the periodic potential $q\in L^2([0,T])$ of $\gamma$, i.e., the regauged complex curvature of $\gamma$ described in Section~\ref{Se:frames}. We also consider the connection form $\alpha = \tfrac12(\lambda\varepsilon + q\varepsilon_+ + \bar{q}\varepsilon_-)$ corresponding to~$q$, and the corresponding extended frame $F(t,\lambda)$. Throughout this proof, we only consider $\lambda\in \bbC$ with $|\lambda-\theta|\leq 1$.

Let $\tilde{q} \in L^2([0,T])$ at first be any periodic potential with $\|\tilde{q}\|_{L^2} \leq \|q\|_{L^2}+1$, let $\tilde{\alpha} = \tfrac12(\lambda\varepsilon + \tilde{q}\varepsilon_+ + \bar{\tilde{q}}\varepsilon_-)$, and let $\tilde{F}$ be the extended frame corresponding to $\tilde{q}$. $\tilde{F}$ is holomorphic in $\lambda$, and we have $\tilde{F}(\cdot,\lambda) \in W^{1,2}([0,T];\mathrm{SL}_2(\bbC)) \subset L^\infty([0,T];\mathrm{SL}_2(\bbC))$ for each $\lambda$, therefore there exists $C_1>0$ (not depending on $\tilde{q}$) so that
\begin{gather}\label{eq:curvedense:F-estim}
\big|\tilde{F}(t,\lambda)\big|, \big|\tilde{F}(t,\lambda)^{-1}\big| \leq C_1
\end{gather}
holds.

Now let $0<\delta\leq 1$ be given. By Corollary~\ref{cor:S3-dense} there exists a finite gap periodic potential $\tilde{q}\in L^2([0,T])$ with $\|q-\tilde{q}\|_{L^2} < \delta$, so that the finite gap curve $\tilde{\gamma}$ corresponding to $\tilde{q}$ is $W^{1,2}$-closed. By definition of the extended frames $F$ and $\tilde{F}$, $D:= F-\tilde{F} \in W^{1,2}\big([0,T];\bbC^{2\times 2}\big)$ is the solution of the ordinary differential equation
\begin{gather}\label{eq:curvedense:D-ode}
\tfrac{{\rm d} \ }{{\rm d}t}D = F\alpha - \tilde{F}\tilde{\alpha} = D\tilde{\alpha} + F(\alpha-\tilde{\alpha})
\end{gather}
with $D(0,\lambda)=0$; by variation of parameters we obtain the explicit formula
\begin{gather}\label{eq:curvedense:D-explicit}
D(t,\lambda) = \int_0^t F(s,\lambda)\bigr(\alpha(s,\lambda)-\tilde{\alpha}(s,\lambda)\bigr)\tilde{F}(s,\lambda)^{-1}\mathrm{d}s \cdot \tilde{F}(t,\lambda) .
\end{gather}
Note that we have $\alpha-\tilde{\alpha} = \tfrac12((q-\tilde{q})\varepsilon_+ + (\bar{q}-\bar{\tilde{q}})\varepsilon_-)$ and therefore
\begin{gather}\label{eq:curvedense:alpha-asymp}
\|\alpha-\tilde{\alpha}\|_{L^2} = \|q-\tilde{q}\|_{L^2} < \delta .
\end{gather}
It follows from equation~\eqref{eq:curvedense:D-explicit} by application of the estimates~\eqref{eq:curvedense:F-estim} (for both $F$ and~$\tilde{F}$) and~\eqref{eq:curvedense:alpha-asymp}
that we have
\begin{gather*}
 |D(t,\lambda)| \leq C_1^3 \|\alpha-\tilde{\alpha}\|_{L_1} \leq C_2 \|\alpha-\tilde{\alpha}\|_{L^2}^{1/2} < C_2 \delta^{1/2}
\end{gather*}
with another constant $C_2>0$. The differential equation~\eqref{eq:curvedense:D-ode} now shows that
\begin{gather*}
 \big\|\tfrac{{\rm d} \ }{{\rm d} t} D(\cdot,\lambda)\big\|_{L^2} < C_3 \delta^{1/2}
\end{gather*}
holds with a constant $C_3>0$ that depends on $\|q\|_{L^2}$, but nothing else. Thus we see that
\begin{gather}
\label{eq:curvedense:D-asymp}
\|D\|_{W^{1,2}} < C_4 \delta^{1/2}
\end{gather}
holds with a constant $C_4>0$.

For $\bbE^3=\bbS^3$ we have $\gamma(t) = F(t,1+\theta)F(t,-1+\theta)^{-1}$ and choose
\begin{gather*} \tilde{\gamma}(t) = \tilde{F}(t,1+\theta)\tilde{F}(t,-1+\theta)^{-1}\end{gather*} according to equation~\eqref{eq:gamma-S3}, and therefore
\begin{gather*}
 \gamma(t)-\tilde{\gamma}(t) = D(t,1+\theta)F(t,-1+\theta)^{-1} \\
 \hphantom{\gamma(t)-\tilde{\gamma}(t) =}{} - \tilde{F}(t,1+\theta)F(t,-1+\theta)^{-1}D(t,-1+\theta)\tilde{F}(t,-1+\theta)^{-1}
\end{gather*}
and
\begin{gather*}
 \gamma'(t) - \tilde{\gamma}'(t) = D(t,1+\theta)\varepsilon F(t,-1+\theta)^{-1} \\
 \hphantom{\gamma'(t) - \tilde{\gamma}'(t) =}{} - \tilde{F}(t,1+\theta)\varepsilon F(t,-1+\theta)^{-1}D(t,-1+\theta)\tilde{F}(t,-1+\theta)^{-1} .
\end{gather*}
It thus follows from the estimates \eqref{eq:curvedense:F-estim} and \eqref{eq:curvedense:D-asymp} that with a constant $C_5>0$,
\begin{gather*}
 \|\gamma-\tilde{\gamma}\|_{W^{1,2}}, \|\gamma'-\tilde{\gamma}'\|_{W^{1,2}} \leq C_1 \delta^{1/2} + C_1^3 \delta^{1/2} \leq C_5 \delta^{1/2}
\end{gather*}
and hence
\begin{gather*}
 \|\gamma-\tilde{\gamma}\|_{W^{2,2}} \leq C_5 \delta^{1/2}
\end{gather*}
holds. Note that $\tilde{\gamma}$ is a periodic finite gap curve with the same period as~$\gamma$, and that $\tilde{\gamma}$ has the same total torsion as $\gamma$ because it was constructed with the same value for~$\theta$.

For $\bbE^3=\bbR^3$ we have $\gamma(t) = 2F'(t,\theta)F(t,\theta)^{-1}$ and choose $\tilde{\gamma}(t) = 2\tilde{F}'(t,\theta)\tilde{F}(t,\theta)^{-1}$ by equation~\eqref{eq:gamma-R3}, where $F'$ and $\tilde{F}'$
again denote the derivative of $F$ and $\tilde{F}$ with respect to $\lambda$, respectively, and therefore
\begin{gather}\label{eq:curvedense:R3-gamma-difference}
\gamma(t)-\tilde{\gamma}(t) = 2D'(t,\theta)F(t,\theta)^{-1} - 2\tilde{F}'(t,\theta)\tilde{F}(t,\theta)^{-1}D(t,\theta)F(t,\theta)^{-1}
\end{gather}
and
\begin{gather}\label{eq:curvedense:R3-gamma'-difference}
\gamma'(t)-\tilde{\gamma}'(t) = 2D'(t,\theta)\varepsilon F(t,\theta)^{-1} - 2\tilde{F}'(t,\theta)\varepsilon\tilde{F}(t,\theta)^{-1}D(t,\theta)F(t,\theta)^{-1} .
\end{gather}
Because $F$, $\tilde{F}$ and $D$ are holomorphic with respect to $\lambda$, we have by Cauchy's inequality
\begin{gather*}
 \|F'(\cdot,\theta)\|_{W^{1,2}} \leq \max_{|\lambda-\theta|=1} \|F(\cdot,\lambda)\|_{W^{1,2}} \overset{\eqref{eq:curvedense:F-estim}}{\leq} C_1 \qquad\text{and likewise}\qquad
 \|\tilde{F}'(\cdot,\theta)\|_{W^{1,2}} \leq C_1
\end{gather*}
and also
\begin{gather*}
 \|D'(\cdot,\theta)\|_{W^{1,2}} \leq \max_{|\lambda-\theta|=1} \|D(\cdot,\lambda)\|_{W^{1,2}} \overset{\eqref{eq:curvedense:D-asymp}}{\leq} C_4 \delta^{1/2} .
\end{gather*}
We now obtain from equations~\eqref{eq:curvedense:R3-gamma-difference}, \eqref{eq:curvedense:R3-gamma'-difference}
\begin{gather*}
 \|\gamma-\tilde{\gamma}\|_{W^{1,2}}, \|\gamma'-\tilde{\gamma}'\|_{W^{1,2}} \leq 2C_1C_4 \delta^{1/2} + 2C_1^3C_4 \delta^{1/2} \leq C_6 \delta^{1/2}
\end{gather*}
with a constant $C_6>0$, and hence
\begin{gather*}
 \|\gamma-\tilde{\gamma}\|_{W^{2,2}} \leq C_6 \delta^{1/2} .
\end{gather*}
Again $\tilde{\gamma}$ is a $T$-periodic, finite gap curve with the same total torsion as $\gamma$. This completes the proof of the theorem.
\end{proof}

\begin{Remark}
 The closed finite gap curves approximating a given closed curve $\gamma$ in Theorem~\ref{T:curvedense} are smooth as a consequence Proposition~\ref{prop:finitegap-periodic}.
 It should be noted, however, that even in the case where the given curve $\gamma$ is in $W^{n,2}$ with $n\geq 3$ (or even smooth),
 the approximation that is claimed by Theorem~\ref{T:curvedense} is only that by the metric of the Sobolev space $W^{2,2}$.

To improve the approximation in this case, one would need to show that the perturbed Fourier map induces local diffeomorphisms $W^{n-2,2}([0,T])_{q,N} \to \ell^2_{n-2} := \big\{ (z_k) \,|\, k^{n-2} z_k \in \ell^2\big\}$ analogous to the local diffeomorphism $\Phi_N$ of Proposition~\ref{prop:Phi-diffeo}, meaning in particular that one would need to prove that the variation of the $(n-2)$-th derivative of the NLS potential $q$ is comparable to the $\ell^2_{n-2}$-measure of the corresponding variation of the perturbed Fourier coefficients. To show such a relationship, it would be necessary to obtain a correspondingly finer control over the asymptotic behaviour of the monodromy and then of the perturbed Fourier coefficients. Such a~control could be obtained by iteratively regauging the connection form $\alpha^q$ from equation~\eqref{eq:alpha_q} to split off terms of a series expansion, however the description of the resulting asymptotic behaviour would be quite complicated.
\end{Remark}

\subsection*{Acknowledgements}
We thank Martin U.~Schmidt for useful discussions and his idea of using perturbed Fourier coefficients. We also thank the referees for their insightful remarks and helpful suggestions that improved the final presentation of the results in this paper. Sebastian~Klein is funded by the Deutsche Forschungsgemeinschaft, Grant 414903103.

\pdfbookmark[1]{References}{ref}
\LastPageEnding

\end{document}